\documentclass[12pt]{imsart}

\RequirePackage[OT1]{fontenc}
\RequirePackage{amsthm,amsmath,amssymb}
\RequirePackage[numbers]{natbib}
\RequirePackage[colorlinks,citecolor=blue,urlcolor=blue]{hyperref}
\usepackage{verbatim}
\usepackage{dsfont}
\usepackage[english]{babel}

\usepackage{epstopdf}
\usepackage{amsfonts} % use AMS math fonts
\usepackage{amsthm}
\usepackage{mathtools}
\usepackage{amsmath}
\usepackage{amssymb}
\usepackage{mathrsfs}
\usepackage{color}
\usepackage{textcomp} % currency symbols
\usepackage{rotating}
\usepackage{graphicx}
\usepackage{caption}
\usepackage{subcaption}
\usepackage{array}
\usepackage{paralist} % package providing inline lists

\arxiv{arXiv:0000.0000}

\startlocaldefs
\numberwithin{equation}{section}
\theoremstyle{plain}

\endlocaldefs

\newcommand{\Cov}{\ensuremath{{\text{Cov}}}} % covariance
\newcommand{\E}{\ensuremath{{\mathbb E}}} % expected value
\newcommand{\R}{\ensuremath{{\mathbb R}}} % real line
\newcommand{\p}{\ensuremath{{\mathbb P}}}

\newcommand{\G}{\mathcal G}
\newcommand{\Var}{\mathrm{Var}}
\newcommand{\N}{\ensuremath{{\mathbb N}}} % natural numbers
\newcommand{\norm}[1]{\left\lVert#1\right\rVert}

\newtheorem{lem}{Lemma}[section]
\newtheorem{Def}{Definition}[section]
\newtheorem{prop}{Proposition}[section]
\newtheorem{theorem}{Theorem}[section]
\newtheorem{rem}{Remark}[section]

\begin{document}

\begin{frontmatter}
\title{Optimal Inference with a Multidimensional Multiscale Statistic}
\runtitle{Inference with a Multiscale Statistic}

\begin{aug}
\author{\fnms{Pratyay} \snm{Datta}\ead[label=e1]{pd2511@columbia.edu}}
\and
\author{\fnms{Bodhisattva} \snm{Sen}\thanksref{t2}\ead[label=e2]{bodhi@stat.columbia.edu}}

\thankstext{t2}{Supported by NSF grant DMS-17-12822 and AST-16-14743.}
%\thankstext{t2}{First supporter of the project}

\runauthor{P.~Datta and B.~Sen}

\affiliation{Columbia University}
\address{Columbia University \\ 1255 Amsterdam Avenue\\
New York, NY 10027\\
\printead{e1}, \; \printead*{e2}}

\end{aug}

\begin{abstract} 
We observe a stochastic process $Y$ on $[0,1]^d$ ($d\geq 1$) satisfying $dY(t)=n^{1/2}f(t)dt$ + $dW(t)$, $t \in [0,1]^d$, where $n \geq 1$ is a given scale parameter (`sample size'), $W$ is the standard Brownian sheet on $[0,1]^d$ and $f \in L_1([0,1]^d)$ is the unknown function of interest. We propose a multivariate multiscale statistic in this setting and prove its almost sure finiteness; this extends the work of~\citet{DC01} who proposed the analogous statistic for $d=1$. We use the proposed multiscale statistic to construct optimal tests for testing $f=0$ versus (i) appropriate H\"{o}lder classes of functions, and (ii) alternatives of the form $f=\mu_n \mathbb{I}_{B_n}$, where $B_n$ is an axis-aligned hyperrectangle in $[0,1]^d$ and $\mu_n \in \R$; $\mu_n$ and $B_n$ unknown.  In the process we generalize Theorem 6.1 of~\citet{DC01} about stochastic processes with sub-Gaussian increments on a pseudometric space, which is of independent interest.
\end{abstract}

\begin{keyword}[class=MSC]
\kwd[Primary ]{62G08}
\kwd{62G86}
\kwd[; Secondary ]{62C20}
\end{keyword}
\begin{keyword}
\kwd{Asymptotic minimax testing}
\kwd{Brownian sheet}
\kwd{kernel estimation}
\kwd{multivariate continuous white noise model}
\kwd{H\"{o}lder classes of functions with unknown smoothness}
\kwd{signal detection}
\end{keyword}

\end{frontmatter}
\section{Introduction}\label{sec: intro}

Let us consider the following continuous multidimensional white noise model:
\begin{equation}\label{eq:Mdl}
{\small  Y(t_1,\ldots,t_d)=\sqrt{n}\int_{0}^{t_1}\ldots\int_{0}^{t_d} f(s_1,\ldots,s_d)\; ds_d \ldots ds_1 + W(t_1,\ldots,t_d),}
\end{equation}
for $(t_1,\ldots,t_d) \in [0,1]^d$ ($d \ge 1$), where $\{Y(t_1,\ldots,t_d): (t_1,\ldots,t_d)\in[0,1]^d \}$ is the observed data, $f \in L_1([0,1]^d)$ is the unknown (regression) function of interest and $W(\cdot)$ is the unobserved $d$-dimensional Brownian sheet (see Definition~\ref{sec:BS}), and $n$ is a known scale parameter. Estimation and inference in this model is closely related to that of nonparametric regression based on sample size $n$.  We work with this white noise model as this formulation is more amiable to rescaling arguments; see e.g.,~\citet{Donoho1992},~\citet{DC01},~\citet{Carter2006}.

In this paper we develop {\it optimal} tests (in an asymptotic minimax sense) based on a newly proposed {\it multidimensional multiscale statistic} (i.e., $d \ge 1$) for testing: 
\begin{itemize}
	\item[(i)] $f=0$ versus a H\"{o}lder class of functions with unknown degree of smoothness; 
	
	\item[(ii)] $f=0$ against alternatives of the form $f=\mu_n \mathbb{I}_{B_n}$, where $B_n$ is an unknown hyperrectangle in $[0,1]^d$ with sides parallel to the coordinate axes (i.e., axis-aligned) and $\mu_n \in \R$ is unknown (for different regimes of $\mu_n$ and $B_n$). 
\end{itemize}
	Scenario (i) arises quite often in nonparametric regression where the goal is to test whether the underlying $f$ is 0 versus $f \ne 0$ with unknown smoothness; see e.g.,~\citet{Lepski2000},~\citet{Horowitz2001},~\citet{Ingster2009multi} and the references therein. Our proposed multiscale statistic, which extends the work of~\citet{DC01} that considered the analogous statistic for $d=1$, leads to rate optimal detection in this problem. Moreover, with the knowledge of the smoothness of the underlying $f$, we construct a {\it asymptotically minimax test} which even attains the exact separation constant (see Section~\ref{Sec 1.2} for formal definitions and related concepts).

Setting (ii) is a prototypical problem in signal detection --- an unknown (constant) signal spread over an unknown hyperrectangular region --- and the goal is to detect the presence of such a signal; see e.g., \citet{Castro2005},~\citet{chan2009},~\citet{Walther10},~\citet{munkchangepoint},~\citet{Butucea2013},~\citet{scanalr},~\citet{Glaz2004},~\citet{Munk2018} for a plethora of examples and applications. 
 
Although several minimax rate optimal tests have been proposed in the literature for this problem (see e.g.,~\citet{Castro2005}, \citet{chan2009}, \citet{Butucea2013} and \citet{Munk2018}), as far as we are aware, our proposed multiscale test is the only test that attains the exact separation constant --- this leads to simultaneous optimal detection of signals both at small and large scales. 

We first motivate and introduce our multiscale statistic below (Section~\ref{secmultiscale}) and briefly describe the asymptotic minimax testing framework and our main optimality results in Section~\ref{Sec 1.2}.

\subsection{Multiscale statistic when $d \ge1$}\label{secmultiscale}
To motivate our multiscale  statistic let us first look at the following testing problem:
\begin{equation} \label{test}
H_0: f=0 \quad \mbox{ versus } \quad H_1: f \neq 0 \in \mathbb{H}_{\beta,L}, 
\end{equation}
where  $\mathbb{H}_{\beta,L}$ is the H\"{o}lder class of function with parameters $\beta >0$ and $L>0$. For $\beta \in (0,1]$ and $L>0$ the H\"{o}lder class $\mathbb{H}_{\beta,L}$ is defined as 
{\small \begin{equation}\label{eq:H_b_L}
\mathbb{H}_{\beta,L}:= \Big \{f\in L_1([0,1]^d):  |f(x) - f(y)| \leq L \norm{x-y}^\beta \mbox{   for all   } x,y \in [0,1]^d \Big\}.
\end{equation} }
For $\beta > 1$ the H\"{o}lder class $\mathbb{H}_{\beta,L}$ is defined similarly; see Definition~\ref{def:holder}.

Our multiscale statistic is based on the idea of {\it kernel averaging}. Suppose that $\psi:\mathbb{R}^d \to \R$ is a measurable function such that (i) $\psi$ is $0$ outside $[-1,1]^d$; (ii) $\psi \in L_2(\R^d)$, i.e., $\int_{\R^d} \psi^2(x) dx <\infty$; (iii) $\psi$ is of bounded (HK)-variation (see Definition~\ref{totalvar}); and (iv) $\int_{\R^d} \psi(x) dx > 0$. We call such a function a {\it kernel}. For any $h:=(h_1,\ldots,h_d) \in (0,1/2]^d$ we define \begin{equation}\label{ah} 
A_h:=\{t\in\mathbb{R}^d : h_i \leq t_i \leq 1-h_i \;\;\mbox{ for  } i =1,\ldots,d\}.
\end{equation} 
For any $t \in A_h$ we define the centered (at $x$) and scaled kernel function $\psi_{t,h}:[0,1]^d \to \R$ as 
{\small \begin{equation}\label{kerest}
\psi_{t,h}(x) :=\psi \left(\frac{x_1-t_1}{h_1},\ldots,\frac{x_d-t_d}{h_d}\right), \quad \mbox{for }\; x=(x_1,\ldots, x_d) \in [0,1]^d.
\end{equation} }
For a fixed $t \in A_h$ we can construct a kernel estimator $\hat{{f_h}}(t)$ of $f(t)$ based on the data process $Y(\cdot)$ as 
\begin{equation*}\label{kerequ-0}
\hat{{f_h}}(t):= \frac{1}{n^{1/2} (\Pi_{i=1}^dh_i) \langle 1,\psi\rangle} \int_{[0,1]^d} \psi_{t,h}(x)dY(x),
\end{equation*}
where for any two functions $g_1,g_2 \in L_2(\R^d)$, we define $\langle g_1,g_2 \rangle:=\int_{\R^d} g_1(x)g_2(x) dx.$
We consider the {\it normalized} version of the above kernel estimator $\hat{{f_h}}(t)$:
\begin{equation}\label{kerest2}
\hat{\Psi}(t,h):=\frac{1}{(\Pi_{i=1}^dh_i)^{1/2}\norm{\psi}}\int_{[0,1]^d} \psi_{t,h}(x) dY(x),
\end{equation}
where $\norm{\psi}^2 := \int_{\R^d} \psi^2(x) dx <\infty$. We can use $\hat{\Psi}(t,h)$ to test $$H_0: f(t)=0 \quad \mbox{ versus } \quad H_1: f(t) \neq 0 $$ where we would reject the null hypothesis for extreme values of $\hat{\Psi}(t,h)$. So, a naive approach to testing \eqref{test} could be to consider $\sup_{ t\in A_h} |\hat{\Psi}(t,h)|$. As this test statistic crucially depends on the choice of the smoothing bandwidth vector $h$, an approach that bypasses the choice of the tuning parameter $h$ and also combines information at various bandwidths would be to consider the test statistic 
\begin{equation}\label{eq:Scan}
\sup_{h>0} \sup_{ t\in A_h} |\hat{\Psi}(t,h)|,
\end{equation} 
where $h>0$ is a short-hand for $h \in (0,\infty)^d$. However, under the null hypothesis $$\sup_{h>0} \sup_{ t\in A_h} |\hat{\Psi}(t,h)|=\infty \qquad \mbox{almost surely (a.s.).}$$  This is because, for a fixed scale  $h$, $\sup_{t \in A_h}|\hat{\Psi}(t,h)| = $ $O_p(\sqrt{2\log(1/(2^d h_1 \cdots h_d))})$; see~\citet{gine02}. Thus, to use the above approach to construct a valid test for~\eqref{test} we need to put the test statistics $\sup_{ t\in A_h} |\hat{\Psi}(t,h)|$ at different scales (i.e., $h$) in the same footing --- this leads to the following definition of the {\it multiscale statistic} in $d$-dimensions:
\begin{equation}\label{eq:TZ}
T(Y,\psi):=\sup_{h\in(0,1/2]^d} \sup_{t \in A_{h}} \frac{\lvert \hat{\Psi}(t,h)\rvert - \Gamma(2^dh_1 \ldots h_d)}{D(2^dh_1 \ldots h_d)}
\end{equation} 
where $\Gamma,D:(0,1] \to (0,\infty)$ are two functions defined as 
\begin{equation*}
\Gamma(r):=(2 \log(1/r))^{1/2}
\end{equation*}
and
\begin{equation*}
D(r):= (\log(e/r))^{-1/2} \log \log (e^e/r).
\end{equation*}  
In Theorem \ref{Thm 1}, a main result in this paper, we show that the above multivariate multiscale statistic $T(Y,\psi)$ is well-defined and finite a.s.~for any kernel function $\psi$, when $f\equiv0$. This result immediately extends the main result of~\citet[Theorem 2.1]{DC01} beyond $d=1$. Although there has been several proposals that extend the definition and the optimality properties of the multiscale statistic  of~\citet{DC01} beyond $d=1$ (see e.g.,~\citet{Munk2018},~\citet{Walther10},~\citet{scanalr}), we believe that our proposed multiscale statistic is the right generalization. Further, the exact form of $T(Y,\psi)$ leads to optimal tests for~\eqref{test} and other alternatives (which the other competing procedures do not necessarily yield; see Remarks~\ref{rem:v} and~\ref{v3} for more details).

To show the finiteness of the proposed multiscale statistic we prove a general result about a stochastic process with sub-Gaussian increments (Theorem~\ref{Thm0}) on a pseudometric space which may be of independent interest. This result has the same conclusion as that of~\citet[Theorem 6.1]{DC01} but assumes a weaker condition on the packing numbers of the pseudometric space on which the stochastic process is defined. This weaker condition on the packing numbers is crucial to the proof of Theorem~\ref{Thm 1}; see Remark~\ref{connection to dumbgen} where we compare our result with the existing result of~\citet[Theorem 6.1]{DC01}. Moreover, Lemma~\ref{lem1} gives a tighter bound on the packing numbers of the pertinent (to our application) pseudometric space, which we believe is also new; see Remarks~\ref{rem:pack} and \ref{rem:v} where we compare our result with some relevant recent work.

\subsection{Optimality of the multiscale statistic}\label{Sec 1.2}

Before we describe our main results let us first introduce the asymptotic minimax hypothesis testing framework. There is an extensive literature on nonparametric testing of the simple hypothesis $\{0\}$. As a staring point we refer the readers to \citet{Ingsterbook}. In the nonparametric setting it is usually assumed that $f$ belongs to a certain class of functions $\mathbb{F}$ and its distance from the null function $f = 0$ is defined by a seminorm $ \norm{\cdot}$. In this setting, given $\alpha \in (0,1)$, the goal is to find a level $\alpha$ test $\phi_n$ (i.e., $\E_0 [\phi_n(Y)] \le \alpha$) such that 
\begin{equation}\label{testinginf}
\inf_{g \in \mathbb{F} : \norm{g} \geq \delta \rho_n}  \E_{g} [\phi_n(Y)]
\end{equation}
is as large as possible for some $\delta >0$ and $\rho_n>0$  where $\rho_n  \to 0$ as $n \to \infty$ ($\rho_n$ is a function of the sample size $n$); in the above notation $\E_g$ denotes expectation under the alternative function $g$. However, it can be shown that given $\mathbb{F}$ and $\norm{\cdot}$, the constants $\delta$ and $\rho_n$ cannot be chosen arbitrarily if one wants to have a statistically meaningful framework (see the survey papers \citet{Ingster19931}, \citet{Ingster19932},   \citet{Ingster19933} for $d=1$ and \citet{Ingster2009multi} for $d>1$). It turns out that if $\delta \rho_n$ is too small then it is not possible to test  the null hypothesis with nontrivial asymptotic power ({i.e., the infimum in~\eqref{testinginf} cannot be strictly larger than $\alpha + o(1)$}). On the other hand if $\delta \rho_n$ is very large many procedures can test $f\equiv 0$ with significant power (i.e., the infimum in~\eqref{testinginf} goes to $1$ as $n \to \infty$). 

The hypothesis testing problem then reduces to: (a) finding the largest possible $\delta \rho_n$ such that no test can have {nontrivial asymptotic power} (i.e., under the alternative $f$ such that $\norm{f}\leq \delta \rho_n  $, the asymptotic  power is less than or equal to the level $\alpha$), and (b) trying to construct test procedures that can detect signals $f$, with $\norm{f} \geq \delta \rho_n$, with considerable power (power going to $1$ as $n \to \infty$). More specifically, $\delta$ and $\rho_n$ are defined such that $\delta \rho_n$ is the largest for which, for all $\epsilon>0$, we have \begin{equation*}
\sup_{\phi_n} \limsup_{n \to \infty}  \inf_{g \in \mathbb{F} : \norm{g} \geq (1-\epsilon)\delta \rho_n}  \E_{g} [\phi_n(Y)] \leq \alpha, 
\end{equation*} 
where the supremum is taken over all sequence of level $\alpha$ tests $\phi_n$. In this case $\rho_n$ is called the {\it minimax rate of testing} and {$\delta$} is called the {\it exact separation constant} (see \citet{Lepski2000}, \citet{Ingster2011sobolev} for more details about minimax testing). On the other hand, we want to find a test $\tilde{\phi}_n$ such that 
\begin{equation*}
\lim_{n \to \infty}\inf_{g \in \mathbb{F} : \norm{g} \geq (1+\epsilon)\delta \rho_n}  \E_{g} [\tilde\phi_n (Y)]=1.
\end{equation*} 
 In such a scenario, $\tilde{\phi}_n$ is called an {\it asymptotically minimax test}.  Here we would also like to point out that if there exists a test $\hat{\phi}_n$ and a constant $\hat{\delta} > \delta$ such that $$\lim_{n \to \infty}\inf_{g \in \mathbb{F} : \norm{g} \geq \hat\delta \rho_n}  \E_{g} [\hat\phi_n(Y)]=1$$ then the test $\hat\phi_n$ is called a {\it rate optimal test}.

In Section~\ref{sec 2.1} we show  that our proposed multiscale statistic yields an asymptotically minimax test for the following scenarios:

\begin{enumerate}
	\item (Optimality for H\"{o}lderian alternatives). Consider testing hypothesis~\eqref{test}. If $$\norm{f}_\infty \geq c_*(1+\epsilon_n)(\log (en)/n)^{\frac{\beta}{2\beta+d}},$$ where $f$ belongs to the H\"{o}lder class $\mathbb{H}_{\beta,L}$ with $ \beta>0$ and $L >0$, $\norm{f}_\infty := \sup_{x \in[0,1]^d} |f(x)|$ denotes the sup-norm of $f$, and $c_*$ is a constant (defined explicitly in Theorem~\ref{multopt}), we show that we can construct a level $\alpha$ test based on the multiscale statistic~\eqref{eq:TZ} that has power converging to 1, as $n \to \infty$, provided $\epsilon_n$ does not go to $0$ too fast (see Theorem~\ref{multopt} for the exact order of $\epsilon_n$). We note that this multiscale statistic would require the knowledge of $\beta$ but not of $L$. $\vspace{0.05in}$
	
	Moreover, we show that if $\norm{f}_\infty \leq c_*(1-\epsilon_n)(\log (en)/n)^{{\beta}/{2\beta+d}}$ no test of level $\alpha \in (0,1)$ can have nontrivial asymptotic power; see  Theorem~\ref{multopt} for the details. {This shows that our proposed multiscale test is asymptotically minimax with rate of testing $\rho_n = (\log(en)/n)^{\beta/(2\beta+d)}$ and exact separation constant $\delta = c_*$. As far as we are aware  this is the first instance of an asymptotically minimax test for the  H\"older class $\mathbb{H}_{\beta,L}$ when $d>1$ (under the  supremum norm). Moreover, if the smoothness $\beta$ of the H\"older class $\mathbb{H}_{\beta,L}$ is unknown (but $\beta \le 1$) then we can still construct a rate optimal test for this problem; see Proposition~\ref{prop:triangle} for the details.}  \vspace{0.07in}

	\item (Optimality for detecting signals at {large/small scales}). Consider testing the hypothesis
	\begin{equation}\label{boxtest}
H_0: f=0 \quad \mbox{ versus } \quad H_1: f=\mu_n \mathbb{I}_{B_n} , 
\end{equation} 
where $\mu_n\neq 0 \in \R$ and $$B_n \equiv B_{\infty}(t^{(n)},h^{(n)}):=\{x \in [0,1]^d : |x_i - t_i^{(n)}| < h_i^{(n)} \mbox{ for all } i=1,\ldots, d\}$$ are unknown, for some $h^{(n)} \in (0,1/2]^d$ and $t^{(n)} \in A_{h^{(n)}}$, and $\mathbb{I}_{B_n}$ denotes the indicator of the hyperrectangle ${B_n}$. First, consider the scenario $\liminf_{n\to \infty} |B_n| > 0 $ where $|B_n|$ denote the Lebesgue measure of $B_n$. Then, if $\lim_{n\to \infty} \sqrt{n} |\mu_n| \to +\infty$, we can construct a level $\alpha$ test based on the multiscale statistic~\eqref{eq:TZ} that has power converging to 1 as $n \to \infty$; see Theorem \ref{boxthm}. Further, we show that, if $\limsup_{n \to \infty} \sqrt{n} |\mu_n| < \infty $, no test of level $\alpha$ can detect the alternative with power going to 1. Thus, the multiscale test is optimal for detecting  signals on large scales. \vspace{0.07in}

On the other hand, let us now consider the case  $\lim_{n \to \infty} |B_n|=0$. If $$|\mu_n|\sqrt{n|B_n|} \geq (1+\epsilon_n) \sqrt{2\log(1/|B_n|)},\quad \mbox{ for all } n,$$ we can construct a test of level $\alpha$, based on the proposed multiscale statistic, that has power converging to 1 as $n \to \infty$, provided $\epsilon_n$ does not go to $0$ too fast (see Theorem \ref{boxthm}). Furthermore, we can show that if $$|\mu_n|\sqrt{n|B_n|} = (1-\epsilon_n) \sqrt{2\log(1/|B_n|)},\quad \mbox{ for all } n,$$ no test can detect the signal reliably with nontrivial power (i.e., for any level $\alpha$ test $\phi_n$ there exists a signal $f_n$ of the above described strength such that $\phi_n$ will fail to detect $f_n$ with asymptotic probability at least $1-\alpha$); see Theorem \ref{boxthm} for the details. This shows that our multiscale test is asymptotically minimax for signals at small scales.
 %Here we would like to point out that the tests constructed for the H\"{o}lderian alternative is different from the test constructed to detect the alternative $f=\mu_n\mathds{B_n}$. We can think of this as 
\end{enumerate}

\subsection{Literature review and connection to existing works}
Our multiscale statistic~\eqref{eq:TZ} can be thought of as a penalized scan statistic, as it is based on the maximum of an ensemble of local test statistics $|\hat{\Psi}(t,h)|$, penalized and properly scaled. Scan-type procedures have  received much attention in the literature over the past few decades. Examples of such procedures can be found in~\citet{Seigmund1995}, \citet{Seigmund2000}, \citet{Naus}, \citet{Kulldorffscan1997}, \citet{haiman2006}, \citet{Jiang2002}, etc. All the above mentioned papers consider $d=1$ and no penalization term (like $\Gamma(\cdot)$ in our case) was used. Asymptotic properties of the scan statistic have been studied expensively. In~\citet{Naus} and~\citet{Pozdnyakov} the authors give asymptotic approximations of the distribution of the scan statistic when $d=1$. For $d=2$, similar results can be found in~\citet{Glaz2004}, \citet{haiman2006},~\citet{WangGlaz2014},  among others. Recently in~\citet{Sharpnack} the authors give exact asymptotics for the scan statistic for any dimension $d$. 

In all of the above papers it is noted that the scan statistic is dominated by small scales; this creates a problem for  detecting large scale signals. One common proposal to fix this problem is to modify the scan statistic so that instead of the maximum over all scales we look at the maximum over scales that are in an appropriate interval containing the true scale of the signal; see e.g.,~\citet{Sharpnack},~\citet{Naus}. In particular, the last two papers show that if the extent of the signal is of a certain order ($\log n$) then this approach leads to power comparable to an oracle. An obvious drawback with the above approach is that we need to have some prior knowledge on which scales the signal(s) may be present. In contrast, our multiscale method does not require any such knowledge. 

 Another approach that has been proposed to optimally detect signals on both large and small scales is to use different critical values (of the scan statistic) to test for signals at different scales separately (see e.g.,~\citet{scanalr},~\citet{Walther10}) and use multiple testing procedures (see~\citet{Hall08} and the references within) to calibrate the method. However, note that a vast majority of the multiple testing literature either assume that the test statistics are independent (which is not the case here) or are too generic and generally quite conservative.

Conceptually, our work is most related to that of~\citet{DC01}, where the authors proposed our multiscale statistic for $d=1$. Thus, our work can be thought of as a generalization of~\citet{DC01} to multidimension ($d>1$).

\subsection{Organization of the paper}
The proposed multiscale statistic is studied in Section~\ref{sec:2}. In Section~\ref{sec 2.1} we construct {\it optimal} tests for: (i) $f=0$ versus H\"{o}lderian alternatives; (ii) $f=0$ versus alternatives of the form $f=\mu_n \mathbb{I}_{B_n}$, where $B_n$ is an axis-aligned hyperrectangle in $[0,1]^d$  and $\mu_n \in \R$ (for different regimes of $\mu_n$ and $B_n$, both unknown). We compare the performance of our multiscale based test with other competing methods in Section~\ref{simulation studies}.  In Section~\ref{future} we discuss some open problems and possible applications/extensions of our work. Section~\ref{proofs of result} gives the proof of Theorem~\ref{Thm 1}. The proofs of the other results are relegated to Appendix~\ref{proofs of result2}.

\section{Multidimensional multiscale statistic}\label{sec:2}
Let us first recall the definition of the multivariate multiscale statistic $T(Y,\psi)$ given in~\eqref{eq:TZ}. The following theorem, our main result in this section, shows that the multiscale statistic $T(Y,\psi)$ is well-defined and finite a.s.~for any (reasonable) kernel function $\psi$; see Section~\ref{sec:Thm 1} for a proof. 
\begin{theorem}\label{Thm 1}
Let $\psi$ be a kernel function. For a positive vector $h:=(h_1,\ldots,h_d)>0$  let $A_{h}$ be as defined in~\eqref{ah}. For $t \in A_h$, let $ \psi_{t,h}(\cdot)$ and $\hat{\Psi}(t,h)$ be as defined in~\eqref{kerest} and~\eqref{kerest2}, respectively. %$$\hat{\Psi}(t,h):=(h_1h_2\ldots h_d)^{-1/2}\norm{\psi}^{-1} \int_{[0,1]^d} \psi_{t,h}(x) dW(x) $$
Consider the statistic $T(W,\psi)$ as defined in~\eqref{eq:TZ}, where $W(\cdot)$ is the Brownian sheet on $[0,1]^d$. Then, almost surely, $T(W,\psi)<\infty$, i.e., $T(W,\psi)$ is a tight random variable.

\end{theorem}

Theorem~\ref{Thm 1} immediately extends the main result of~\citet[Theorem 2.1]{DC01} beyond $d=1$. The proof of the above theorem crucially relies on the following two results. We first introduce some notation.\begin{Def}[Packing number] \label{packing}
For any pseudometric space $(\mathscr{F},\rho)$ and $\epsilon>0$, the packing number $N(\epsilon,\mathscr{F})$ is defined as the supremum of the number of elements in $\mathscr{F}^\prime$ where $\mathscr{F}^\prime \subseteq \mathscr{F}$ and for all $a \ne b \in \mathscr{F}^{\prime}$ we have $\rho(a,b)> \epsilon.$
\end{Def}

We will prove Theorem~\ref{Thm 1} as a consequence of the following more general result about stochastic processes with sub-Gaussian increments on some pseudometric space (see Section~\ref{pf:Thm0} for its proof). 

\begin{theorem}\label{Thm0}
Let $X$ be a stochastic process on a pseudometric space  $(\mathscr{F},\rho)$ with continuous sample paths. Suppose that the following three conditions hold:
\begin{itemize}
\item[(a)] There is a function $\sigma : \mathscr{F} \to (0,1]$ and a constant $K \geq 1$ such that $$\p \big(X(a) > \sigma (a) \eta \big) \leq K \exp(-\eta^2/2) \qquad \forall \, \eta > 0,\; \forall \, a\in \mathscr{F}.$$ 
Moreover, $\sigma^2(b) \leq \sigma^2(a) + \rho^2(a,b),\;\; \forall \;a,b \in \mathscr{F}$.

\item[(b)] For some constants $L,M \geq 1$, $$\p\big(|X(a)-X(b)| > \rho(a,b)\eta\big) \leq L \exp(-\eta^2/M)\quad \forall \, \eta > 0, \; \forall\, a,b\in \mathscr{F} .$$

\item[(c)] For some constants $A,B,V, p > 0$,
$${\small N((\delta u)^{1/2}, \{a \in \mathscr{F}: \sigma^2(a) \leq \delta\})\leq A u^{-B} \delta^{-V} (\log(e/\delta))^{p} \;\; \forall \, u, \delta \in (0,1].}$$ 
Then the random variable 
\begin{equation}\label{eq:S(X)}
S(X):= \sup_{a \in \mathscr{F}} \frac{X^2(a)/\sigma^2(a) - 2V\log(1/\sigma^2(a))}{\log \log (e^e/\sigma^2(a))}
\end{equation} is finite almost surely. More precisely, $\p(S(X) > r) \le \xi(r)$ for some function $\xi:\R_+ \to \R$ depending only on the constants $K,L,M,A,B,p,V$ such that $\lim_{r \to \infty} \xi(r) = 0$.
\end{itemize}
\end{theorem}
\begin{rem}[Connection to~\citet{DC01}]\label{connection to dumbgen}
A similar result to Theorem~\ref{Thm0} above appears in~\citet[Theorem 6.1]{DC01}. However note that there is a subtle and important difference: The bound on the packing number in (c) of Theorem~\ref{Thm0} involves the additional logarithmic factor $(\log(e/\delta))^{p}$ which is not present in~\citet[Theorem 6.1]{DC01}. In fact, we show that even with this additional logarithmic factor, the random variable $S(X)$, defined in~\eqref{eq:S(X)}, involves the same penalization term $2V\log(1/\sigma^2(a))$ as in~\citet[Theorem 6.1]{DC01}. Hence, we can think of Theorem~\ref{Thm0} as an generalization of~\citet[Theorem 6.1]{DC01}. 
\end{rem}
To apply Theorem~\ref{Thm0} to prove Theorem~\ref{Thm 1} we need to define a suitable pseudometric space $(\mathscr{F},\rho)$ and a stochastic process, and verify that conditions (a)-(c) in Theorem~\ref{Thm0} hold. In that vein, let us define the following set 
$$\mathscr{F} := \left\{(t,h) \in \mathbb{R}^d \times (0,1/2]^d:  h_i\leq t_i\leq 1-h_i, \mbox{ for all } i=1,2,\ldots,d \right\}$$
with the following pseudometric $$\rho^2((t,h),(t^\prime,h^\prime)):=|B_\infty(t,h) \bigtriangleup B_\infty(t^\prime,h^\prime)|, \quad \quad \mbox{for  } (t,h), (t^\prime,h^\prime) \in \mathscr{F},$$
where $B_\infty(t,h) :=\Pi_{i=1}^d(t_i-h_i,t_i+h_i)$, $A \bigtriangleup B := (A \cap B^c)\cup (A^c \cap B)$ denotes the symmetric difference of the sets $A$ and $B$, and $|A|$ denotes the Lebesgue measure of the set $A$. Also, define $$\sigma^2(t,h):= |B_\infty(t,h)|=2^d \Pi_{i=1}^d h_i, \qquad \mbox{for  } (t,h) \in \mathscr{F}.$$ 
The following important result shows that indeed for the above defined pseudometric space $(\mathscr{F},\rho)$ condition (c) of Theorem~\ref{Thm 1} holds.
\begin{lem}\label{lem1}
Let $\mathscr{F},\rho(\cdot,\cdot)$  and $\sigma(\cdot)$ be as described above. Then 
\begin{equation*}
N\left((u\delta)^{1/2},\{(t,h) \in \mathscr{F}: \sigma^2(t,h)\leq \delta\}\right) \leq K u^{-2d}\delta^{-1} (\log(e/\delta))^{d-1} \;\; \forall \, u,\delta \in(0,1], 
\end{equation*}
for some constant $K$ depending only on $d$.
\end{lem}
\begin{rem}\label{rem:pack}
Here we would like to point out that Lemma~\ref{lem1} shows that condition (c) of Theorem~\ref{Thm0} holds with $B=2d$, $p=d-1$ and most importantly for  $V=1$, which was also the case when $d=1$ (as shown in~\citet{DC01}).
\end{rem}
\begin{rem}\label{rem:v}
Compare the numerator of our multiscale statistic~\eqref{eq:TZ} with the multiscale statistic proposed in~\citet[Equation (6)]{Munk2018}. In \citet{Munk2018} the authors propose a penalization term $\Gamma_V(2^dh_1\ldots h_d)$ where $\Gamma_V:(0,1]\to (0,\infty)$ is defined as $$\Gamma_V(r):=(2V \log(1/r))^{1/2}$$
instead of the penalization $\Gamma(2^dh_1\ldots h_d)$ as in~\eqref{eq:TZ}. Further, in that paper the authors recommend the choice of $V=(2d-1+\epsilon)$ for any $\epsilon>0$; see~\citet[Lemma 5.1]{Munk2018}. Thus, Theorem~\ref{Thm 1} and Lemma~\ref{lem1}, improve on the existing results in the literature. Our penalization term $\Gamma(\cdot)$ results in optimal detection properties for testing \eqref{test} and \eqref{boxtest} which cannot be achieved if the penalization term $\Gamma_V(\cdot)$, for $V >1$, is used.
\end{rem}
It is well-known that we should choose the constant $V$ in the penalization term $\Gamma_V $ as small as possible (see e.g.,~\citet[Section 1.1]{Munk2018}) for optimal testing. In our proposed multiscale statistic we take $V=1$. The following proposition shows that indeed $V=1$ is the smallest possible permissible value; see Section~\ref{sec:Prop_V} for a proof.
\begin{prop}\label{v<1}
Suppose $V<1$. Let $\Gamma_V$ and $\mathscr{F}$ be as defined above. Then we have $$\sup_{(t,h)\in \mathscr{F}} |\hat{\Psi}(t,h)|-\Gamma_V(2^dh_1\ldots h_d)=\infty\quad \mbox{a.s.}$$
Thus, $\sup_{(t,h)\in \mathscr{F}} \frac{|\hat{\Psi}(t,h)|-\Gamma_V(2^dh_1\ldots h_d)}{D(2^dh_1\ldots h_d)}=\infty\quad \mbox{a.s.}$
\end{prop}

\section{Optimality of the multiscale statistic in testing problems}\label{sec 2.1}
In this section we prove that we can construct tests based on the multiscale statistic that are optimal for testing~\eqref{test} and~\eqref{boxtest}. For both the testing problems we can define a multiscale test based on kernel $\psi $ as follows:
Let $$\kappa_{\alpha,\psi}=\inf\{c \in \R: \p(T(W,\psi)> c) \leq \alpha\},$$
where $W$ is the standard Brownian sheet on $[0,1]^d$. For notational simplicity we would denote $\kappa_{\alpha,\psi}$ by $\kappa_\alpha$ from now on.

For testing~\eqref{test} and~\eqref{boxtest} a test of level $\alpha$ can be defined as follows: $$\mbox{Reject $H_0\quad\quad$ if and only if $\qquad T(Y,\psi)>\kappa_{\alpha}$}.$$ Let us call this testing procedure the multiscale test. Although any kernel $\psi$ can be used to construct the above test, in Sections~\ref{sec 3.1} and~\ref{sec 3.2} we show that specific choices of the kernel function $\psi$ leads to asymptotically minimax tests. 

\subsection{Optimality against H\"older classes of functions}\label{sec 3.1}
Let us recall the definition of the H\"{o}lder class of functions $\mathbb{H}_{\beta,L}$, for $\beta \in (0,1]$ and $L>0$, as in~\eqref{eq:H_b_L}; see Definition~\ref{def:holder} for the formal definition of $\mathbb{H}_{\beta,L}$ for any $\beta >0$. Let $\psi_\beta:\R^d \to \R$, for $0<\beta< \infty$, be the unique solution of the following optimization problem:
\begin{equation}\label{minimize} \mbox{Minimize} \norm{\psi} \mbox{ over all } \psi \in \mathbb{H}_{\beta,1} \mbox{ with } \psi(0) \geq 1.\end{equation} Elementary calculations show that for $0< \beta \leq 1$, we have $$\psi_\beta(x)=(1-\norm{x}^\beta) \mathbb{I}(\norm{x} \leq 1);$$ see Section~\ref{lem:Ele} for a proof. For $ \beta>1$, $\psi_\beta$ can be calculated numerically. We consider the kernel $\psi_\beta$, for $\beta>0$, described above and state our first optimality result for testing \eqref{test}; see Section~\ref{multiopt} for a proof.

\begin{theorem}\label{multopt}
Let $T_\beta \equiv T(Y,\psi_\beta)$ be the multiscale  statistic defined in~\eqref{eq:TZ} with kernel $\psi_{\beta}$, for $0<\beta<\infty$. Define $$\rho_n :=\left(\frac{\log n}{n}\right)^{\frac{\beta}{2\beta+d}}$$ and $$c_* \equiv c_*(\beta,L) :=\left(\frac{2dL^{d/\beta}}{(2\beta+d)\norm{\psi_\beta}^2}\right)^{\frac{\beta}{2\beta+d}}.$$
Then, for arbitrary $\epsilon_n > 0 $ with $\epsilon_n \to 0$ and $\epsilon_n \sqrt{\log n }\to \infty$ as $n \to \infty$, the following hold:
\begin{itemize}
	\item[(a)] For any arbitrary sequence of tests $\phi_n$ with level $\alpha$ for testing~\eqref{test}, we have 
$$\limsup_{n \to \infty} \inf_{g\in \mathbb{H}_{\beta,L}: \norm{g}_{\infty}= (1-\epsilon_n)c_*\rho_n} \E_g [\phi_n(Y)] \leq \alpha;$$

	\item[(b)] for $J_n :=[(c_*\rho_n/L)^{1/\beta},1-(c_*\rho_n/L)^{1/\beta}]^d$, we have 
$$\lim_{n\to \infty} \inf_{g\in \mathbb{H}_{\beta,L}: \norm{g}_{J_n,\infty}\geq (1+\epsilon_n)c_*\rho_n} \p_g(T_\beta > \kappa_\alpha)=1$$
where $\norm{g}_{J_n,\infty}:= \sup_{t \in J_n}|g(t)|$.
\end{itemize}
\end{theorem}
The above result generalizes~\citet[Theorem 2.2]{DC01} beyond $d=1$. Theorem~\ref{multopt} can be interpreted as follows: (a) for every test $\phi_n$ there exists a function with supremum norm $(1-\epsilon_n)c_*\rho_n$ which cannot be detected with nontrivial asymptotic power; whereas  (b) when we restrict to functions with signal strengths (i.e., supremum norm in the interior of $[0,1]^d$) just a bit larger than the above threshold, our proposed multiscale test is able to detect every such function with asymptotic power 1. In this sense our proposed test is optimal in detecting departures from the zero function for H\"{o}lder classes $\mathbb{H}_{\beta,L}$. We note here that to calculate $T_\beta$  we need the knowledge of $\beta$ but we do not need to know $L$. 

If $\beta$ is unknown, but is less than or equal to 1, we can use $T_1$ as a test statistic for testing~\eqref{test}. Although the resulting test is not asymptotically minimax, the test is still rate optimal. The following result formalizes this; see Section~\ref{proof:proptriangle} for its proof.

\begin{prop}\label{prop:triangle}
Consider testing~\eqref{test} where $\beta\leq 1$ is unknown. Let us recall the definition of $\psi_1$ in~\eqref{minimize}. Let $T_1 \equiv T(Y,\psi_1)$ be the multiscale  statistic defined in~\eqref{eq:TZ} with kernel $\psi_{1}$. Define $$\rho_n :=\left(\frac{\log n}{n}\right)^{\frac{\beta}{2\beta+d}}$$ and let $M$ be any constant such that $M >\left(\frac{2dL^{d/\beta}\norm{\psi_1}^2}{(2\beta+d)\langle \psi_1,\psi_\beta \rangle^2}\right)^{\frac{\beta}{2\beta+d}}.$
Let $J_n :=[(M\rho_n/L)^{1/\beta},1-(M\rho_n/L)^{1/\beta}]^d$. Then we have 
$$\lim_{n\to \infty} \inf_{g\in \mathbb{H}_{\beta,L}: \norm{g}_{J_n,\infty}\geq M\rho_n} \p_g(T> \kappa_\alpha)=1$$
%where $\norm{g}_{J_n,\infty}:= \sup_{t \in J_n}|g(t)|$.
where $\kappa_\alpha$ is the $(1-\alpha)$ quantile of the multiscale statistic $T(Y,\psi_1)$ under the null hypothesis.
\end{prop}
\begin{rem}
Instead of using the test statistic $T_\beta$ if we use the test statistic
\begin{equation}\label{eq:T^*}T_\beta^\star:=\sup_{h\in(0,1/2]^d} \sup_{t \in A_h} \big[\lvert \hat{\Psi}(t,h)\rvert - \Gamma(2^dh_1 \ldots h_d)\big]
\end{equation} with the kernel $\psi_\beta$, then the same conclusions as that of Theorem~\ref{multopt} and Proposition~\ref{prop:triangle} would hold. Thus the multiscale statistic $T_\beta^\star$ is also optimal against H\"{o}lderian alternatives. 
\end{rem}
\begin{rem}\label{rem:Gamma_V}
Note that in~\citet{Munk2018} the authors propose a multiscale statistic like $T_\beta^\star$, with a slightly different penalization term 
\begin{equation}\label{eq:Gamma_V} 
\Gamma_V: r \mapsto (2V \log(1/r))^{1/2}
\end{equation} 
instead of $\Gamma(\cdot)$. A close inspection of our proof of Theorem~\ref{multopt} reveals that for such a statistic, only signals with $\norm{g}_{J_n,\infty}\geq \sqrt{V}(1+\epsilon_n)c_*\rho_n$ will be detected with power converging to 1. This shows how a proper penalization (as in our multiscale statistic) can lead  to the testing procedure attaining the exact separation constant for testing~\eqref{test}.
\end{rem}
\subsection{Optimality against axis-aligned hyperrectangular signals}\label{sec 3.2}

In Theorem~\ref{multopt} we proved the optimality of the multiscale test when the supremum norm of the signal is large. A natural question that arises next is: ``What if the signal is not peaked but distributed evenly on some subset of $[0,1]^d$?". To answer this question  we look at the testing problem~\eqref{boxtest}, and establish below the optimality of our multiscale test in this setting (see Section~\ref{multiopt} for a proof of Theorem~\ref{boxthm}). Note that when $d=1$ similar optimality results are known for the multiscale statistic; see~\citet[Theorem 2.6]{munkchangepoint} and~\citet{scanalr}. For $h=(h_1,\ldots,h_d) \in (0,1/2]^d$, let us first define $$\mathscr{B}_{h}:=\{B \subseteq [0,1]^d: B=\Pi_{i=1}^d[t_i - h_i,t_i+h_i] \mbox{ for some } t=(t_1,\ldots,t_d) \in A_h \}.$$
\begin{theorem}\label{boxthm}
Let $T \equiv T(Y,\psi_0)$ where $\psi_0= \mathbb{I}_{[-1,1]^d}$. Let $f_n=\mu_n\mathbb{I}_{B_n}$ where $B_n$ is an axis-aligned hyperrectangle and let $|B_n|$ denote the Lebesgue measure of the set $B_n$. Then we have the following results: 
\begin{itemize}
\item[(a)] Suppose that $\liminf_{n \to \infty }|B_n| > 0$. Let $\phi_n$ be any test of level $\alpha\in(0,1)$ for \eqref{boxtest}. Then, for any $f_n = \mu_n\mathbb{I}_{B_n}$  such that  $\limsup_n |\mu_n|\sqrt{n |B_n|}< \infty$, we have $$\limsup_{n \to \infty}  \E_{f_n} [\phi_n(Y)] <1.$$ Moreover, for the proposed multiscale test based on $T$, we have $$\lim_{n \to \infty} \inf_{f_n:\lim |\mu_n|\sqrt{n |B_n|}=\infty} \p_{f_n}(T>\kappa_\alpha)=1.$$

\item[(b)]Now let us look at the case $\lim_{n\to \infty} |B_n|=0$. Let $h_n=(h_{1,n},\ldots,h_{d,n}) \in (0,1/2]^d$ be any sequence of points such that $\lim_{n \to \infty} \Pi_{i=1}^d h_{i,n} \to 0$. Let $$\G_n^{-}:=\{f_n=\mu_n\mathbb{I}_{B_n}:|\mu_n|\sqrt{n|B_n|} = (1-\epsilon_n)\sqrt{2\log(1/|B_n|)} , B_n \in \mathscr{B}_{h_n}\}$$ with $\epsilon_n \to 0$ and $\epsilon_n \sqrt{2\log(1/|B_n|)} \to \infty$. (Here we have omitted the dependence of $h_n$ in the notation $\G_n^{-}$). If $\phi_n$ be any test of level $\alpha\in(0,1)$ for \eqref{boxtest}  then we have $$ \limsup_{n \to \infty} \inf_{f_n \in \G_n^{-}} \E_{f_n} [\phi_n(Y)] \leq \alpha.$$
Moreover, let $$\G_n^{+}:=\{f_n=\mu_n\mathbb{I}_{B_n}:|\mu_n|\sqrt{n|B_n|} \geq (1+\epsilon_n)\sqrt{2\log(1/|B_n|)},B_n \in \mathscr{B}_{h_n}\}.$$ Then for our multiscale test we have $$\lim_{n \to \infty} \inf_{f_n \in \G_n^{+}} \p_{f_n}(T>\kappa_\alpha)=1.$$
\end{itemize}

\end{theorem}

\begin{rem}
If we use the test statistic $T^\star$, as defined in~\eqref{eq:T^*} (with the kernel $\psi_0$), instead of $T$ in Theorem~\ref{boxthm}, the optimality results  described in the theorem still hold.
\end{rem}
Our first result in Theorem~\ref{boxthm} shows that as long as $\liminf_{n \to \infty }|B_n| > 0$, for any test to have power converging to $1$ we need to have $\lim |\mu_n|\sqrt{n |B_n|}=\infty$, in which case our multiscale test achieves asymptotic power 1. Thus our multiscale test is optimal for detecting large scale signals. The next result can be interpreted as follows: (i) For signals with small spatial extent (i.e., $\lim_{n \to \infty} |B_n|=0$) if the signal strength is too small ($|\mu_n|\sqrt{n|B_n|} \le (1-\epsilon_n)\sqrt{2\log(1/|B_n|)}) $ no test can detect the signal reliably with nontrivial probability (i.e., for every test $\phi_n$ there exist a signal such that $\phi_n$ will fail to detect it with probability $1-\alpha+o(1)$); (ii) on the other hand, if the signal strength is a bit larger than the threshold (i.e., the exact separation constant) described above our multiscale test will detect the signal with asymptotic power 1. This shows that  our multiscale test achieves optimal detection for signals with small spatial footprint. We would like to emphasize here that by using the same exact test (using the same kernel $\psi_0$) we are able to optimally detect both large and small scale signals.

\begin{rem}\label{v3} 
As we mentioned in Remark~\ref{rem:Gamma_V} if we used $\Gamma_V(\cdot)$ (see~\eqref{eq:Gamma_V}), for $V>1$, instead of $\Gamma(\cdot)$, in defining the multiscale statistic then we would only be able to detect signals (when $|B_n| \to 0$) if $|\mu_n|\sqrt{n|B_n|} \geq \sqrt{V}(1+\epsilon_n) \sqrt{2\log(1/|B_n|)}$ which is not the exact separation constant as mentioned in Theorem~\ref{boxthm}. This agains illustrates the importance of choosing the right penalization term $\Gamma(2^dh_1\ldots h_d)$ in defining the multiscale statistic.
\end{rem}
\begin{rem}
Here we would like to point out that proofs for the minimax lower bound that have been derived for the two scenarios in Theorems ~\ref{multopt} and~\ref{boxthm} follows the standard techniques that have been used in  \citet{Ingster19931}, \citet{Ingster19932}, \citet{Ingster19933},~\citet{Lepski2000}, \citet{DC01}, \citet{Ingster2009multi} etc.  
\end{rem}

\subsubsection{Comparison with the scan and average likelihood ratio statistics when $d=1$}
When $d=1$ there exists an extensive literature on the optimal detection threshold for signals of the form $f_n=\mu_n\mathbb{I}_{B_n}$, where now $B_n \subseteq [0,1]$ is an interval. In~\citet{scanalr} the authors compare the performance of the scan statistic (i.e., the statistic~\eqref{eq:Scan} in the discrete setup with $\psi = \mathbb{I}_{[-1,1]}$) and the average likelihood ratio (ALR) statistic (which is the discrete analogue of $\int_0^{1/2} \int_{h}^{1-h} \exp[ |\hat{\Psi}(t,h)|^2/2] dt \, dh$); see Section~\ref{simulation studies} for a description and comparison of the two competing methods with our multiscale test when $d=2$. 

When $\liminf_{n \to \infty} |B_n|>0$ the scan statistic can only detect the signal, with asymptotic power 1, when $|\mu_n|\sqrt{n} \geq (1+\epsilon_n)\sqrt{2\log n}$, whereas the ALR statistic (and the proposed multiscale statistic) can detect the signal whenever we have $|\mu_n|\sqrt{n} \to \infty$ (which is a less stringent condition). Note that $|\mu_n|\sqrt{n} \to \infty$ is also required for any test to detect the signal with asymptotic power $1$. This shows that the scan statistic is not optimal for detecting large scale signals.

On the other hand if $\lim_{n \to \infty} |B_n|=0$, the scan statistic can detect the signal if $|\mu_n|\sqrt{n|B_n|} \geq (1+\epsilon_n) \sqrt{2\log n}$ whereas the ALR statistic can detect the signal when $|\mu_n|\sqrt{n|B_n|} \geq \sqrt{2}(1+\epsilon_n) \sqrt{2\log(1/|B_n|)}$.  The optimal detection threshold in this scenario is $ |\mu_n|\sqrt{n|B_n|} \geq (1+\epsilon_n) \sqrt{2\log(1/|B_n|)}$, which is attained by the multiscale statistic. Thus that scan statistic is optimal in detecting signals only when $|B_n|=O(1/n)$. The ALR statistic requires the signal to be at least $\sqrt{2}$ times the (detectable) threshold. This shows that neither the standard scan or the ALR is able to achieve the optimal threshold for detecting small scale signals. 
 
 \citet[Theorem 2.6]{munkchangepoint} shows the optimality of the multiscale statistic (which is a modification of the scan statistic) in detecting signals in both cases when $d=1$. In \citet{Rivera13} and \citet{scanalr} the authors propose a condensed ALR statistic which, much like the multiscale statistic, is able to attain the optimal threshold for detection in both regimes of $B_n$. As far as we are aware the condensed ALR statistic has not been extended beyond $d=1$ and therefore whether it achieves the optimal threshold for $d>1$ is not known. 
 In summary, Theorem~\ref{boxthm} shows that our multidimension multiscale test is asymptotically minimax even when $d>1$.

\section{Simulation studies}\label{simulation studies}
\begin{figure}
\includegraphics[scale=.50]{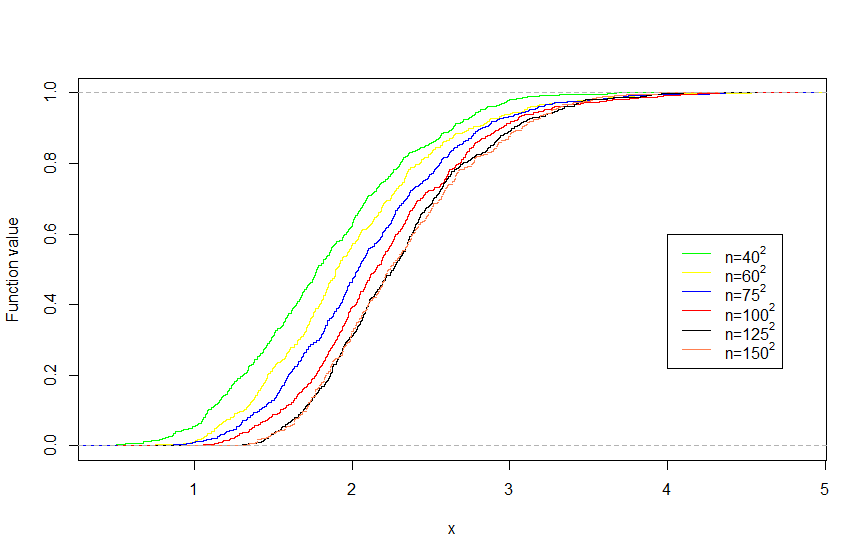}\caption{The empirical distribution functions of the multiscale statistic for different values of $n$.}\label{multi dist}
\end{figure}

\begin{table}
\begin{tabular}{| c | c || c| c|}
\hline 
\multicolumn{4}{|c|}{Critical values}\\
\hline
$m$ & 95\% quantile & $m$ & 95\% quantile \\
\hline
25 & 3.02 & 75 & 3.27  \\
40 & 3.12  &  100 & 3.31 \\
50 & 3.18 & 125 & 3.32 \\
60 & 3.22  & 150  & $~3.30^\star$ \\
\hline
\end{tabular}
\caption{Critical values $\kappa_{0.05}$ for different $n = m^2$.}\label{tab:Crit}
{\it  $^\star$Note that 0.95 quantiles necessarily  increase as $n$ increases. But in our simulations the 0.95 quantile for $n=150^2$ turned out to be slightly less than that of $n=125^2$ due to sampling variability.}
\end{table}  
In this section we demonstrate the performance of the multiscale testing procedure described in Section~\ref{sec 2.1} and compare it with other competing methods through simulation studies. For computational tractability, we replace the continuous white noise model~\eqref{eq:Mdl} with a discrete one and consider the case $d=2$. More specifically, we consider data on the $m\times m$ grid  $S_n = \{(i/m,j/m): 1\leq i,j\leq m\}$ (here $n =m^2$), where the model is 
\begin{equation*}
Y\left(\frac{i}{m},\frac{j}{m}\right)=f\left(\frac{i}{m},\frac{j}{m}\right) + \epsilon\left(\frac{i}{m},\frac{j}{m}\right), \qquad \mbox{for }\; i, j = 1,\ldots, m,
\end{equation*}
with $\epsilon(i/m,j/m)$'s being i.i.d.~standard normal random variables.  In our simulation experiments we vary our bandwidth parameter $h=(h_1,h_2)$ in the $m \times m$ grid $S_n$. For the simulations we have used the kernel function $\psi =\mathbb{I}_{[-1,1]^d}$. In Table~\ref{tab:Crit} we give the empirical 0.95-quantile of the multiscale statistic $T(W,\psi)$ (see~\eqref{eq:TZ}) for different values of $n$; the computation of the empirical quantiles were based on 3000 replications. Observe that the empirical quantiles seem to stabilize as $m$ increases beyond 100. Figure~\ref{multi dist} shows the empirical distribution function estimates, based on 3000 replications, of the multiscale statistic for different values of $n$. 

In Tables~\ref{Tab: com1} and~\ref{Tab: com2} we compare the powers of the multiscale test, a test based on a scan-statistic, and the  ALR test (see \citet{scanalr} for the details). Formally, we consider testing~\eqref{boxtest} against alternatives of the form $H_1: f=\mu_n \mathbb{I}_{B_n}$, for both small and large scale signals ($B_n$). We briefly describe the above two competing procedures. Let $\mathscr{B}$ be the set of all axis-aligned rectangles on $[0,1]^2$ with corner points of the form $({i}/{m},{j}/{m})$, for $i,j \in \{1,\ldots, m\}$. For every $B \in \mathscr{B}$ define $$\hat{\Psi}(B) := \frac{1}{\sqrt{|B|}} \sum_{(i/m,j/m)\in B} Y\left(\frac{i}{m},\frac{j}{m}\right).$$ Note that $\hat{\Psi}(\cdot)$ is the discrete analogue of the normalized kernel estimator as defined in~\eqref{kerest2}. The scan test statistic (see \citet[Chapter 5]{Scanbook}) for this problem is defined as $$ M_n:= \max_{B \in \mathscr{B}} |\hat{\Psi}(B)|. $$ The ALR test statistic (see \citet{chan2009}) is defined as $$A_n:= \frac{1}{{m \choose 2}^2} \sum_{B \in \mathscr{B}} \exp(\hat{\Psi}(B)^2/2).$$ The scan test (ALR test) rejects the null hypothesis if the observed $M_n$ ($A_n$) exceeds the 0.95-quantile for $M_n$ ($A_n$) under the null hypothesis. In Tables~\ref{Tab: com1} and \ref{Tab: com2} we compare the performance of the three procedures. Here $\mu$ denotes the signal strength, and $k/m$ denotes the length of each side of the square signal $B_n$ (here $m=40$ and $100$ for the two cases). The power of the tests were  calculated using 1000 replications.
\begin{table}

\begin{tabular}{| c | c | c| c| c|}
\hline 
\multicolumn{4}{|c|}{$k=1$}\\
\hline 
$\mu$ & Scan & Multiscale & ALR\\
\hline
3.5 & 0.23 &  0.08  & 0.07 \\
4.0 & 0.34  &  0.13  & 0.08 \\
4.5 & 0.50  &  0.18   & 0.08 \\
5.0& 0.71 & 0.30  & 0.08 \\
5.5 & 0.86 & 0.53  & 0.09 \\
\hline

\end{tabular} \quad 
\begin{tabular}{| c | c | c| c| c|}
\hline 
\multicolumn{4}{|c|}{$k=4$}\\
\hline 
$\mu$ & Scan & Multiscale  & ALR\\
\hline
1.00 & 0.22  &  0.14  & 0.11 \\
1.20 & 0.43 &  0.31  &  0.30 \\
1.35 & 0.60 & 0.48  & 0.44 \\
1.50 & 0.74 & 0.55    & 0.52 \\
1.65 & 0.86  & 0.72  & 0.61  \\
\hline

\end{tabular}

\vspace{.3cm}

\begin{tabular}{| c | c | c| c| c|}
\hline 
\multicolumn{4}{|c|}{$k=18$}\\
\hline 
$\mu$ & Scan & Multiscale & ALR\\
\hline
0.20 & 0.15  & 0.21  & 0.19  \\
0.30 &  0.49 & 0.68  &  0.67 \\
0.35 & 0.65 & 0.80  &  0.82 \\
0.40 & 0.80  & 0.90  &  0.89 \\
\hline

\end{tabular} \quad
\begin{tabular}{| c | c | c| c| c|}
\hline 
\multicolumn{4}{|c|}{$k=40$}\\
\hline 
$\mu$ & Scan & Multiscale & ALR\\
\hline
0.040 & 0.15  & 0.32  & 0.31 \\
0.043 & 0.30  & 0.56  &  0.54 \\
0.047 & 0.45  & 0.78  & 0.78 \\
0.050 & 0.68  & 0.94  & 0.95  \\
\hline

\end{tabular}

\caption{Power of the scan, the multiscale and the ALR tests for $m=40$ (i.e., $n=40^2$) as $\mu$ changes.}
\label{Tab: com1}
\end{table}

We make the following observations. For both the cases ($m=40$ and $100$)  when the signal is at the smallest scale, e.g., $k=1$, the scan statistic outperforms everything else. However, when $m=100$, even in relatively small scales, e.g., $k=8$ (i.e., about $0.6\%$ of the observations contain the signal) our multiscale test starts to outperform the scan test. Note that in this setting (small scales) the ALR performs the worst. As the spatial extent of the signal increases, our multiscale procedure and the ALR procedure starts performing favorably whereas the performance of the scan statistics deteriorates. Thus, the simulation experiments corroborates our theoretical findings.

\begin{table}

\begin{tabular}{| c | c | c| c| c|}
\hline 
\multicolumn{4}{|c|}{$k=1$}\\
\hline 
$\mu$ & Scan & Multiscale & ALR\\
\hline
4.5 & 0.34 &  0.11  & 0.06 \\
5.0 & 0.52  &  0.28  & 0.06 \\
5.5 & 0.75  &  0.43   & 0.09 \\
6.0 & 0.95 & 0.61  & 0.13 \\
\hline

\end{tabular} \quad 
\begin{tabular}{| c | c | c| c| c|}
\hline 
\multicolumn{4}{|c|}{$k=8$}\\
\hline 
$\mu$ & Scan & Multiscale  & ALR\\
\hline
0.25 & 0.08  &  0.17  & 0.07 \\
0.30 & 0.35 &  0.46  &  0.13 \\
0.35 & 0.60 & 0.72  & 0.22 \\
0.40 & 0.82 & 0.96    & 0.50 \\
\hline

\end{tabular}

\vspace{.3cm}

\begin{tabular}{| c | c | c| c| c|}
\hline 
\multicolumn{4}{|c|}{$k=30$}\\
\hline 
$\mu$ & Scan & Multiscale & ALR\\
\hline
0.040 & 0.07  & 0.22  & 0.22  \\
0.050 &  0.17 & 0.42  &  0.45\\
0.055 & 0.42 & 0.74  &  0.75 \\
0.060 & 0.58  & 0.93  &  0.96 \\
\hline

\end{tabular} \quad
\begin{tabular}{| c | c | c| c| c|}
\hline 
\multicolumn{4}{|c|}{$k=100$}\\
\hline 
$\mu$ & Scan & Multiscale & ALR\\
\hline
0.014 & 0.08  & 0.42  & 0.42 \\
0.018 & 0.17  & 0.62  &  0.63 \\
0.020 & 0.22  & 0.84  & 0.86 \\
0.025 & 0.45  & 0.96  & 0.95  \\
\hline

\end{tabular}
\caption{Power of the scan, the multiscale and the ALR tests for $m=100$ (i.e., $n=100^2$) as $\mu$ changes.}
\label{Tab: com2}
\end{table}

 \section{Discussion}\label{future}
In this paper we have proposed a multidimensional multiscale statistic in the  continuous white noise model and used this statistic to construct asymptotically minimax tests for testing $f=0$ against (i) H\"{o}lder classes of functions; and (ii) alternatives of the form $f=\mu_n \mathbb{I}_{B_n}$, where $B_n$ is an unknown axis-aligned hyperrectangle in $[0,1]^d$ and $\mu_n \in \R$ is unknown. However, there are many open questions in this area. We briefly delineate a few of them below and in the process describe some important papers in related areas of research.

We have shown that for the H\"older class $\mathbb{H}_{\beta,L}$, if the smoothness parameter $\beta$ is known, we can construct an asymptotically minimax test. However, if $\beta$ is unknown (and $\beta \le 1$) we can only construct a rate optimal test. A natural question that arises is whether a test can be constructed that is asymptotically minimax (for the H\"{o}lder class of functions with the supremum norm) without the knowledge of the smoothness parameter $\beta$ (and $L>0$); see~\citet[Section 1.3]{Ji2017}. Another interesting question would be to try to extend our results to other smoothness classes like Sobolev/Besov classes; in \citet{Ingster2011sobolev} the authors gave the minimax rate of testing for Sobolov class, but no test was proposed that achieves the exact separation constant.

Note that we have shown that our multiscale test is asymptotically minimax for detecting the presence of a signal on an axis-aligned hyperrectangle in $[0,1]^d$. One obvious extension of our work would be to correctly identify the hyperrectangle on which the signal is present. Further, we could go beyond hyperrectangles and try to identify signals that are present on some other geometric structures $A \subset [0,1]^d$ (i.e., $f=\mu \mathbb{I}_A$ where $A$ is not necessarily an axis-aligned hyperrectangle). Examples of such geometric structures could be: $(i)$ $A$ is an hyperrectangle which is not necessarily axis-aligned, $(ii)$ $A$ is a $d$-dimensional ellipsoid, $(iii)$ $A=\bigcup_{i=1}^k A_i$ where each $A_i \subseteq [0,1]^d$ is an (axis-aligned) hyperrectangle, etc.~\citet{munkchangepoint} and the references therein investigated the problem of finding change points in $d=1$ which can be thought of as detection of multiple intervals. In~\citet{Castro2005} the authors use the scan statistic to detect regions in $\mathbb{R}^d$ where the underlying function is non-zero.~\citet{Arias2010} considers the problem of finding a cluster of signals (not necessarily rectangular) in a network using the scan statistic. Although the method they propose achieves the optimal boundary for detection, it requires the knowledge of whether the signal shape is ``thick" or ``thin". For  hyperrectangles this refers to whether or not the minimum side length is of order $\log n /n$ or not. We believe that the multiscale statistic, with proper modifications, can be used to find asymptotically minimax/rate optimal tests in such problems.

In our white noise model~\eqref{eq:Mdl} we assume  that the distribution of the response variables is (homogeneous and independent) Gaussian. Similar questions about signal detection can be asked when the response is non-Gaussian; see e.g.,~\citet{Munk2018},~\citet{ChanWalther2015}, \citet{Rivera13}, \citet{Walther10} etc. In~\citet{Munk17heterochange} the authors looked at the problem of detecting change points under heterogeneous variance of the response variable (when $d=1$).~\citet{Rohde2008} looked at this problem where the error distribution is known to be symmetric (when $d=1$). \citet{Walther10} studied a similar problem where the response variable is binary. A multiscale approach could be used to tackle such problems as well. 

Several interesting applications of the multiscale approach exist when $d=1$ (following the seminal paper of~\citet{DC01}): In~\citet{Dumbgen2008} the authors propose a multiscale test statistic to make inference about a probability density on the real line  given i.i.d.~observations;~\citet{Hieber2013} use multiscale methods to make inference in a deconvolution problem;~\citet{Rivera13} use multiscale methods to detect a jump in the intensity of a Poisson process, etc. We believe that our extension beyond $d=1$ will also lead to several interesting multidimensional applications.

\section*{Acknowledgements}
The authors would like to thank Lutz D\"{u}mbgen and Sumit Mukherjee for several helpful discussions.
\section{Proofs of our main results}\label{proofs of result}
\subsection{Some useful concepts}\label{useful concepts}
In this subsection we formally define some technical concepts that we use in this paper.
\begin{Def}[Brownian sheet]\label{sec:BS}
By a $d$-dimensional Brownian sheet we mean a mean-zero Gaussian process $\{W(t): t \in [0,1]^d\}$ with covariance
\begin{equation*}
\Cov(W(t_1,\ldots,t_d),W(s_1,\ldots,s_d))=\Pi_{i=1}^d \min({t_i,s_i}),
\end{equation*} 
for $(t_1,\ldots,t_d), (s_1,\ldots,s_d)\in [0,1]^d$. The Brownian sheet is the $d$-dimensional counterpart of the standard Brownian motion; see e.g.,~\citet{WZ01},~\citet[Chapter 5]{Bsheetbook} for detailed properties of the Brownian sheet.
\end{Def}

In the following we give some useful properties of a Brownian sheet $W(\cdot)$.
\begin{itemize}
\item If $g \in L_2([0,1]^d) $ then $\int g dW:=\int_{[0,1]^d} g(t) dW(t) \sim N(0,\norm{g}^2).$

\item If $g_1,g_2 \in L_2([0,1]^d) $ then $\Cov \left(\int g_1 dW,\int g_2 dW\right)= \int_{[0,1]^d} g_1(t)g_2(t) dt. $

\item { {\it Cameron-Martin-Girsanov Theorem for Brownian sheet:} Let us state the simplest version of the Cameron-Martin-Girsanov Theorem that we will use in this paper (see~\citet[Chapter 3]{Protterbook} for detailed discussion about change of measure and the result). $\vspace{0.05in}$

Assume $f \in L_1([0,1]^d)$ and let $\{{W(t):t\in[0,1]^d}\}$ be a standard Brownian sheet. Let $\Omega$ be the set of all real-valued continuous functions defined on $[0,1]^d$. Let $P$ denote the measure on $\Omega$ induced by the Brownian sheet  $\{{W(t):t\in[0,1]^d}\}$ and let $Q$ denote the measure induced by $\{Y(t): t \in [0,1]^d\} $ where $Y(t)$ is defined as in \eqref{eq:Mdl}. 
Then $Q$ is absolutely continuous with respect to $P$ and the Radon-Nikodym  derivative is given by  $$\frac{dQ}{dP}(Y)=\exp\left(\sqrt{n}\int fdW - \frac{n}{2}\norm{f}^2 \right).$$ This, in turn, implies that for any measurable function $\phi$ we have  $$ \E_Q\left(\phi(Y)\right) = \E_P \left(\phi(Y)\frac{dQ}{dP}(Y)  \right).$$}
\end{itemize}

Let us now define the H\"{o}lder class of functions $\mathbb{H}_{\beta,L}$, for $\beta>0$ and $L>0$.
\begin{Def}\label{def:holder}
Fix $\beta>0$ and $L>0$. Let $\lfloor \beta \rfloor$ be the largest integer which is strictly less than $\beta$ and for $k=(k_1,k_2,\ldots ,k_d)\in \mathbb{N}^d$ set $\norm{k}_1 := \sum_{i=1}^d k_i$. The H\"{o}lder class $\mathbb{H}_{\beta,L}$ on $[0,1]^d$  is the set of all functions $f:[0,1]^d \to \R$ having all partial derivatives of order $\lfloor \beta \rfloor$ on $[0,1]^d$ such that $$\sum_{0\leq \norm{k}_1 \leq \lfloor \beta \rfloor} \sup_{x \in [0,1]^d} \left| \frac{\partial^{\norm{k}_1} f(x) }{\partial x_1^{k_1} \ldots \partial x_d^{k_d}} \right| \leq L$$ and $$\sum_{\norm{k}_1 = \lfloor \beta \rfloor} \left| \frac{\partial^{\norm{k}_1} f(y) }{\partial x_1^{k_1} \ldots \partial x_d^{k_d}} -  \frac{\partial^{\norm{k}_1} f(z) }{\partial x_1^{k_1} \ldots \partial x_d^{k_d}} \right| \leq L \norm{y-z}^{\beta -\lfloor \beta \rfloor } \quad \forall  \, y,z \in [0,1]^d. $$
\end{Def}

\begin{rem} One of the most important properties of $\mathbb{H}_{\beta,L}$ that we will use is the following: If $f  \in \mathbb{H}_{\beta,1}$ then, for any $h=(h_1,\ldots,h_d)>0$ and $t\in A_h$,  $$g(x_1,\ldots, x_d) := L \min(h)^\beta f\left(\frac{x_1-t_1}{h_1}, \ldots ,\frac{x_d-t_d}{h_d} \right) \in \mathbb{H}_{\beta,L}$$
where $\min(h):= \min_{i=1,\ldots,d} h_i$.

\end{rem}
 
\begin{Def}[Hardy-Krause variation]\label{totalvar} 
The notion of bounded variation for a function $f: \R^d \to \R$, where $d \ge 2$, is more involved than when $d=1$. In fact there is no unique notion of bounded variation for a function when $d \ge 2$. Below we describe the  notion of Hardy and Krause variation as given in~\citet{Variationpaper}, which suffices for our purpose. 

Let $f: [-1,1]^d \to \R$ be a measurable function. Let $a=(a_1,\ldots,a_d)$ and $b=(b_1,\ldots,b_d)$ be elements of $[-1,1]^d$ such that $a<b$ (coordinate-wise). We introduce the d-dimensional difference operator $\Delta^{(d)}$ which assigns to the axis-aligned box $A:=[a,b]$ a d-dimensional quasi-volume
$$\Delta^{(d)}(f;A)=\sum_{j_1=0}^1\cdots \sum_{j_d=0}^1 (-1)^{j_1+\cdots+j_d}f(b_1+j_1(a_1-b_1),\ldots,b_d+j_d(a_d-b_d)).$$
Let $m_1, \ldots, m_d \in \N$. For $s=1,\ldots,d$, let $-1=:x_0^{(s)}<x_1^{(s)}<\cdots < x_{m_s}^{(s)}:=1$ be a partition of $[-1,1]$ and let $\mathsf{P}$ be a partition of $[-1,1]^d$ which is given by
$$\mathsf{P}:=\left\{[x_{l_1}^{(1)},x_{l_1+1}^{(1)}] \times \cdots \times [x_{l_d}^{(d)},x_{l_d+1}^{(d)}]: \; l_s=0,1,\ldots,m_s-1, \;\mbox{for } s=1,\ldots,d \right\}.$$
Then the variation of $f$ on $[-1,1]^d$ in the sense of {\it Vitali} is given by
$$V^{(d)}(f;[-1,1]^d):=\sup_{\mathsf{P}} \sum_{A \in \mathsf{P}}|\Delta^{(d)}(f;A)|$$
where the supremum is extended over all partitions of $[-1,1]^d$ into axis-parallel boxes generated by $d$ one-dimensional partitions of $[-1,1]$. For $1 \leq s \leq d$ and $1\leq i_1 < \ldots < i_s \leq d$, let $V^{(s)}(f;i_1,\ldots,i_s;[-1,1]^d)$ denote the $s$-dimensional variation
in the sense of Vitali of the restriction of $f$ to the face $$U_d^{(i_1,\ldots,i_s)}=\left\{ (x_1,\ldots,x_d) \in [-1,1]^d : x_j = 1 \mbox{ for all } j \neq i_1,\ldots, i_s \right\}$$ of $[-1,1]^d$. Then the variation of $f$ on $[-1,1]^d$ in the sense of Hardy and Krause anchored at 1, abbreviated by HK-variation, is given by $$TV(f):= \sum_{i=1}^d \sum_{1 \leq s \leq d} V^{(s)}(f;i_1,\ldots,i_s;[-1,1]^d).$$
We say a function $f$  has   bounded HK-variation if $TV(f) < \infty.$
\end{Def}
The main property of a bounded HK-variation function that we will need in this paper is stated below.
\begin{rem}
If $f$ is a right continuous function on $[-1,1]^d$ which has bounded HK-variation then there exists a unique signed Borel measure $\nu$ on $[-1,1]^d$ for which $$f(x)=\nu([-1,x]), \quad x \in [-1,1]^d;$$see \citet{Variationpaper} for details.
\end{rem}

\subsection{Proof of Theorem~\ref{Thm0}}\label{pf:Thm0}
In the following proofs $K$ would be used to denote a generic constant whose value would change from line to line. 

For every $v > 0$, we define $$\Gamma(X,v):=\sup_{a,b \in \mathscr{F},\rho(a,b)\leq v} |X(a)-X(b)|.$$ For simplicity we divide the proof in three steps.\newline

\noindent \textbf{Step 1:} In this step we will prove that \begin{equation}\label{tailsup}
\p\big(\Gamma(X,v) > \eta \big) \leq K \exp \left( -\frac{\eta^2}{Kv^2\log(e/v)} \right) \quad \forall \, \eta >0 \mbox{ and } v \in (0,1], 
\end{equation} where $K>0$ is a positive constant not depending on $v$. We will prove the above result by introducing the notion of Orlicz norm. Let $\lambda:\mathbb{R}_+ \to \mathbb{R}$ be a nondecreasing convex function with $\lambda(0)=0$. For any random variable $X $ the  Orlicz norm $\norm{X}_{\lambda}$ is defined as $$\norm{X}_{\lambda}= \inf\bigg\{C>0: \E\lambda \left(\frac{|X|}{C}\right)\leq 1\bigg\}.$$ The Orlicz norm is of interest to us as any Orlicz norm easily yields a bound on the tail probability of a random variable i.e., $\p(|X|>x) \leq [\lambda(x/\norm{X}_{\lambda})]^{-1}, \mbox{ for all } x \in \R.$
Let us define $\lambda(x) :=\exp(x^2)-1$, $x >0$.
Hence,
\begin{equation}\label{eq:Orlicz}
\p \big(|X|>x\big) \leq \min \bigg\{1,\frac{1}{\exp(x^2/\norm{X}^2_{\lambda})-1}\bigg\} \leq 2\times \exp(-x^2/\norm{X}^2_{\lambda}).
\end{equation} Hence, it is enough to bound the Orlicz norm of $\Gamma(X,v)$. A bound on the Orlicz norm of $\Gamma(X,v)$ can be shown by appealing to  \citet[Theorem 2.2.4]{VW01} which we state below. 

\begin{lem}\label{lem:Cov_no}
 Let $\lambda:\R_+\to \R$ be a convex, nondecreasing, non-zero function with $\lambda(0)=0$  and for some constant $c>0$, $\limsup_{x,y \to \infty} \frac{\lambda(x)\lambda(y)}{\lambda(cxy)} < \infty $. Let $\{X_a, a \in \mathscr{F}\}$ be a separable stochastic process with $$\norm{X_a - X_b}_{\lambda} \leq C \rho(a,b) \mbox{ for all } a,b \in \mathscr{F} $$ for some pseudometric $\rho$ on $ \mathscr{F}$ and constant C. Then for any $\zeta , v > 0$,
 $$ \norm{\Gamma(X,v)}_{\lambda} \leq K \left[\int_{0}^{\zeta} \lambda^{-1}(N(\epsilon,\mathscr{F})) d\epsilon + v \lambda^{-1}(N^2(\zeta,\mathscr{F}))\right] $$ for some constant $K$ depending only on $\lambda$ and $C$. 
 \end{lem} 
 
We apply the above lemma with $\lambda(x) :=\exp(x^2)-1$ (i.e., $\lambda^{-1}(y)=\sqrt{\log(1+y)}$). Note that condition (b) of Theorem~\ref{Thm0} directly implies that $\norm{X_a - X_b}_{\lambda} \leq C \rho(a,b)$ by an application of~\citet[Lemma 2.2.1]{VW01}. 

By taking $\delta=1,\epsilon=u^{1/2}$, condition (c) of Theorem~\ref{Thm0} yields $N(\epsilon,\mathscr{F}) \leq A\epsilon^{-2B}$. Thus, Lemma~\ref{lem:Cov_no} gives (with $\zeta = v$)
 $$ \norm{\Gamma(X,v)}_{\lambda} \leq K \left[ \int_{0}^{v} \sqrt{\log(1+A\epsilon^{-2B})} d\epsilon + v \sqrt{\log(1+A^2v^{-4B})} \right].$$
The expression on the right side of the above display can be easily shown to be less than or equal to $ K v \sqrt{\log(e/v)}$ for some constant $K$. This result along with an application of~\eqref{eq:Orlicz} with $\Gamma(X,v)$ instead of $X$ imply 
 $$\p\big(\Gamma(X,v) > \eta \big) \leq K \exp\left(-\frac{\eta^2}{Kv^2\log(e/v)}\right) \qquad \mbox{for all } \eta >0, \; 0<v\leq1,$$ for some constant $K$.\newline
 
\noindent \textbf{Step 2:} Let us define $\mathscr{F(\delta)}:= \{a \in \mathscr{F}: \delta/2 < \sigma^2(a) \leq \delta\}$, for $\delta \in (0,1]$, and 
\begin{equation}\label{eq:Pi-delta}\Pi(\delta):=\p\left(\frac{X^2(a)}{\sigma^2(a)} > 2V\log(\frac{1}{\delta})+S\log\log(\frac{e^e}{\delta}) \mbox{ for some } a\in \mathscr{F}(\delta)\right)
\end{equation}
for $S \geq 4p+1$. In this step we will prove that $$\Pi(\delta) \leq K \exp((K-S/K)\log\log(e^e/\delta))$$ for some constant $K$.

Fix $u < 1/2$. Let $\mathscr{F}(\delta,u)$ be a $\sqrt{u \delta}$-packing set of $\mathscr{F(\delta})$. By our assumption the cardinality of $\mathscr{F}(\delta,u)$ is less than or equal to $Au^{-B}\delta^{-V} (\log({e/\delta}))^p$. Fix $a \in \mathscr{F}(\delta)$. From the definition of $\mathscr{F}(\delta,u)$ we can associate $\hat{a} \in \mathscr{F}(\delta,u)$ (corresponding to $a \in \mathscr{F}(\delta)$) such that $ \rho^2(a,\hat{a}) \leq u\delta$. Using assumption (a) of Theorem~\ref{Thm0} we have \begin{equation}\label{equ:ineq}
\sigma^2(a) \geq \sigma^2(\hat{a}) - u\delta \geq \sigma^2(\hat{a})(1-2u)\end{equation} where the last inequality follows from the fact that $\hat{a} \in \mathscr{F(\delta})$ (thus $\sigma^2(\hat{a}) > \delta/2$).  

We want to study the event 
\begin{equation}\label{event}
 \frac{X^2(a)}{\sigma^2(a)} > r 
\end{equation} for some $r>0$. Obviously, 
for any $\lambda \in (0,1)$, either (i) $|X(a) - X(\hat{a})|^2 > \lambda^2 X^2(a)$ or (ii) $|X(a) - X(\hat{a})|^2 \le \lambda^2 X^2(a)$ (which, in particular implies $|X(\hat{a})| \ge (1 - \lambda) |X(a)|$). The above two cases reduce to: 
\begin{equation}\label{eventcase1}
 \Gamma(X,(u\delta)^{1/2})^2 \geq |X(a) - X(\hat{a})|^2 > \lambda^2 X^2(a) \geq \lambda^2 r \sigma^2(a) \geq\lambda^2 r \frac{\delta}{2} 
\end{equation} 
(here the first inequality follows from the definition of $\Gamma(X,(u\delta)^{1/2})$ and the third inequality follows from condition \eqref{event}), and
\begin{equation}\label{eventcase2}
X^2(\hat{a}) \geq (1-\lambda)^2 X^2(a) \geq  (1-\lambda)^2 r \sigma^2(a) \geq (1-\lambda)^2 r (1-2u) \sigma^2(\hat{a}) 
\end{equation} (here the second inequality follows from \eqref{event} and last inequality follows from \eqref{equ:ineq}). Therefore, for any $r >0$,
\begin{eqnarray*} 
\Pi_r(\delta) &:=  &\p\left(\frac{X^2(a)}{\sigma^2(a)} > r \mbox{  for some } a \in \mathscr{F(\delta)}\right)\\ 
& \leq  & \p\left(\Gamma(X,(u\delta)^{1/2})^2 > \lambda^2 \delta r/2\right) \\ 
& & \qquad \qquad \qquad  + \sum_{\hat{a} \in \mathscr{F}(\delta,u)} \p\left(X^2(\hat{a})/\sigma^2(\hat{a}) > (1-\lambda)^2 r (1-2u)\right)
\end{eqnarray*} 
where we have used the fact that if $X^2(a)/\sigma^2(a) > r$ for some $a \in \mathscr{F}$, then either~\eqref{eventcase1} holds or~\eqref{eventcase2} is satisfied for some $\hat{a} \in \mathscr{F}(\delta,u)$. The first term on the right side of the above display can be bounded by appealing to~\eqref{tailsup} with $\eta=\sqrt{\lambda^2\delta r /2}$ and $v=\sqrt{u\delta}$ and the second term can be bounded by using conditions (a) and (c) of Theorem~\ref{Thm0}. Hence we get
\begin{align} 
\Pi_r(\delta) & \leq   K \exp\left(- \frac{\lambda^2 \delta r/2}{Ku\delta \log(e/\sqrt{u\delta})} \right) \nonumber \\ 
& \qquad \qquad \qquad+ A u^{-B}\delta^{-V} \big(\log(\frac{e}{\delta})\big)^p \exp\left(-\frac{(1-\lambda)^2 r (1-2u)}{2}\right) \nonumber \\
%& & \Big[ \mbox{The first term follows from \eqref{tailsup} by appropriately choosing }  v,\eta \mbox{ and } \\ \nonumber 
%& &   \mbox{ the second term follows from condition  } (a) \mbox{ and } (c) \mbox{ of Theorem } \ref{Thm0} \Big] \\ \nonumber
& \leq K \Big[\exp\left(-\frac{\lambda^2 r}{K u \log(e/(u\delta))}\right) \nonumber \\
& \qquad  + \exp\big(B \log(1/u)+V\log(1/\delta)+p\log\log(e/\delta)+ur- (1/2-\lambda)r \big)\Big]. \label{eq1} 
\end{align}

Fix $S \geq 8p+1$ and set $$r:= 2V \log(1/\delta) + S \log\log\big(\frac{e^e}{\delta}\big)$$ and  $$\lambda:=\frac{1}{r} \Big((S/4)\log\log(e^e/\delta) - p\log \log (e/\delta)\Big).$$ Observe that $r >1$ and $0<\lambda< 1/4$. Moreover, we have $$(1/2-\lambda)r= V\log(1/\delta)+p\log\log(e/\delta)+ (S/4)\log\log(e^e/\delta).$$
Putting these values in~\eqref{eq1} gives us 
\begin{eqnarray}
\Pi(\delta) \equiv \Pi_r(\delta)
& \leq &  K \Bigg[\exp\left(-\frac{(S-4p)^2 (\log\log(e^e/\delta))^2}{Kur\log(e/(u\delta))}\right) \nonumber \\
& &\quad \quad \quad+ \exp\Big(B\log(1/u)+ ur- (S/4)\log\log(e^e/\delta)\Big) \Bigg] \qquad \label{eq2}
\end{eqnarray}
where we have used the fact that $\lambda^2 r^2 = ((S/4)\log\log(e^e/\delta) - p\log \log (e/\delta))^2 \ge (S-4p)^2 (\log\log(e^e/\delta))^2/16$.
Now, let us pick $$u:= \frac{S}{8r \log(e/\delta)} < \frac{1}{2}.$$
Then we have $\frac{1}{u} \leq K \log^2(e/\delta)$ for some constant $K$. Let us consider the two terms on the right side of~\eqref{eq2} separately. For the first term, using $ur = S[ \log(e/\delta)]^{-1}/8$, and that $\frac{1}{u} \leq K \log^2(e/\delta)$,  we have
{\small \begin{eqnarray*}
\frac{(S-4p)^2 (\log\log(e^e/\delta))^2}{Kur\log(e/(u\delta))}  & = &\frac{8(S-8p + 16p^2/S) (\log\log(e^e/\delta))^2\log(e/\delta)}{K\big(\log(e/\delta) + \log (u^{-1})\big)}\\
%& \ge & \frac{(S-8p) (\log\log(e^e/\delta))^2\log(e/\delta)}{K\big(\log(e/\delta) + \log K + 2\log \log (e/\delta)\big)}\\
& \ge &(S-8p)(\log\log(e^e/\delta))\Big(\frac{(\log\log(e^e/\delta))\log(e/\delta)}{K\big(\log(e/\delta) + \log K + 2\log \log (e/\delta)\big)}\Big)\\
& \geq  & (1/K')(S-8p)(\log\log(e^e/\delta)).
\end{eqnarray*}}
Here the last inequality follows from the following fact: As $$\tau(\delta) := \frac{(\log\log(e^e/\delta))\log(e/\delta)}{K\big(\log(e/\delta) + \log K + 2\log \log (e/\delta)\big)} \to \infty, \qquad \mbox{ as }\delta \to 0,$$ we can find a lower bound $K'>0$ such that $\tau(\delta) \ge 1/K'$ for all $\delta \in (0,1]$. 

For the second term on the right side of~\eqref{eq2} we have 
\begin{eqnarray*}
&& B\log(1/u)+ur -(S/4)\log\log(e^e/\delta) \\
& \leq &  B\log K+2B\log\log(e/\delta)+S/8-(S/4)\log\log(e^e/\delta)\\
& \leq & B\log K+2B\log\log(e/\delta) - (S/8)\log\log(e^e/\delta)\\
& \leq & B\log K+(2B-S/8)\log\log(e^e/\delta).\end{eqnarray*} 
Thus, both the terms on the right side of~\eqref{eq2} have the form $K \exp[(C - S/K') \log\log(e^e/\delta)]$ for some constants $K,C,K' >0$. Putting these values in~\eqref{eq2} gives us, for suitable constant $K >0$, we get $$\Pi(\delta)\leq K \exp\left((K-S/K)\log\log(e^e/\delta)\right).$$

\noindent \textbf{Step 3:} In this step we will prove that as $S \to \infty$ $$\p\left(X^2(a)/\sigma^2(a) > 2V\log(1/\sigma^2(a)) + S \log\log\big(\frac{e^e}{\sigma^2(a)}\big)\mbox{ for some } a \in \mathscr{F} \right) \to 0.$$
First let us define $$\tilde{\Pi}(\delta) :=\p\left(X^2(a)/\sigma^2(a) > 2V\log(1/\sigma^2(a)) + S \log\log\big(\frac{e^e}{\sigma^2(a)}\big)\mbox{ for some } a \in \mathscr{F}(\delta) \right).$$ 
Comparing with~\eqref{eq:Pi-delta} we can see that for any $\delta \in (0,1]$, $$\tilde{\Pi}(\delta) \leq \Pi(\delta) $$ as: If $a \in \mathscr{F}(\delta)$ then  $\sigma^2(a) \leq \delta$ and $x \longmapsto 2V\log(1/x)+ S \log \log (e^e/x)$ is a decreasing function of $x$. Hence, we have 
$$\tilde{\Pi}(\delta)\leq K \exp\left((K-S/K)\log\log(e^e/\delta)\right).$$
Therefor, as $\mathscr{F}=\bigcup_{l\geq 0} \mathscr{F}(2^{-l})$, 
\begin{align*}
\p \Big(X^2(a)/\sigma^2(a) &  >  2V \log(1/\sigma^2(a)) + S \log\log(\frac{e^e}{\sigma^2(a)})\mbox{ for some } a \in \mathscr{F} \Big)\\
& \leq \sum_{l=0}^\infty \tilde{\Pi}(2^{-l})\\
& \leq K  \sum_{l=0}^\infty  \exp((K-S/K)\log\log(e^e2^l))\\
&= K \sum_{l=0}^\infty (e+l\log 2)^{-(S/K-K)} \quad
\to 0 \quad \mbox{   as   } S\to \infty.
\end{align*}
This proves that $S(X) :=\sup_{a \in \mathscr{F}} \frac{X^2(a)/\sigma^2(a) - 2V\log(1/\sigma^2(a))}{\log \log (e^e/\sigma^2(a))} < \infty$ a.s.\qed %\newline

\subsection{Proof of Lemma~\ref{lem1}}\label{pf:lem1}
First let us define the following sets:
\begin{eqnarray*}
 \mathscr{F}_{\delta,(l_1,\ldots ,l_d)} & :=  &\big\{(t,h) \in \mathscr{F}: \delta/2 < \sigma^2(t,h)\leq  \delta,  \; 2^{l_i-1}  < \frac{h_i}{\delta^{1/d}} \leq 2^{l_i}, \;\forall \; i=1,\ldots ,d \big\} \\
& & \hspace{2.5in}  \mbox{ for some }(l_1,\ldots,l_d) \in \mathbb{Z}^d, \\
 \mathscr{F}(\delta) & := & \big\{(t,h) \in \mathscr{F}: \delta/2 < \sigma^2(t,h)\leq  \delta \big\}.
\end{eqnarray*}
We note  that $\mathscr{F}_{\delta,(l_1,\ldots ,l_d)}$ is empty unless we have 
$$ \mathrm{(i)} \quad l_i \leq (1/d) \log_2(1/\delta) \qquad \mbox{for all } i=1,\ldots ,d;$$ (this restriction is a consequence of the fact that $h_i \leq 1/2$) and $$ \mathrm{(ii)}  \quad -(d+1) < \sum_{i=1}^d l_i \leq 0$$ (this restriction is a consequence of the fact that $\delta/2 < \sigma^2(t,h) \leq \delta$).\newline

\noindent{\bf Step 1:}
First, we will show that for any $(l_1,\ldots , l_d) \in \mathbb{Z}^d$, and $\delta,u \in (0,1]$,
\begin{equation}\label{eq3}
N\left((u\delta)^{1/2},\mathscr{F}_{\delta,(l_1,\ldots,l_d)}\right)  \leq   K u^{-2d}\delta^{-1}.
\end{equation}
Let $\mathscr{F}^{\prime}$ be a subset of $\mathscr{F}_{\delta,(l_1,\ldots,l_d)}$ such that for any two elements  $(t,h),(t^\prime,h^\prime) \in \mathscr{F}^\prime$ we have \begin{equation}\label{packingcond} \rho^2((t,h),(t^\prime,h^\prime)) > u\delta.\end{equation} Our aim is to show that $$|\mathscr{F}^\prime| \leq Ku^{-2d}\delta^{-1} ,$$ for some constant $K$ independent of $(l_1,\ldots,l_d)$, $u$ and $\delta$. If $\mathscr{F}_{\delta,(l_1,\cdots,l_d)} $ is empty then the assertion is trivial.  So assume that $\mathscr{F}_{\delta,(l_1,\cdots,l_d)} $ is non-empty which imposes bounds on the $l_i$'s as shown above.

Let us define the following partition of $[0,1]^d$ into disjoint hyperrectangles:
\begin{eqnarray*} 
R:=\Big\{M_{(i_1,\ldots, i_d)}\cap [0,1]^d: M_{(i_1,\ldots, i_d)} :=\Pi_{k=1}^d \Big((i_k-1)\frac{u\delta^{\frac{1}{d}}2^{l_k}}{c}, i_k\frac{u\delta^{\frac{1}{d}}2^{l_k}}{c} \Big], \\
 \qquad \qquad \: 1 \le i_k \le \lceil cu^{-1}\delta^{-\frac{1}{d}}2^{-l_k} \rceil \Big\}
\end{eqnarray*} 
where we take $c :=d4^d$.
We would like to point out that in the above definition when $i_k=1$, for any $k=1,\ldots,d$, by  $\big((i_k-1)c^{-1}u\delta^{1/d}2^{l_k},$ $ i_kc^{-1}u\delta^{1/d}2^{l_k}\big]$ we mean the closed interval $\big[0, c^{-1}u\delta^{1/d}2^{l_k}\big]$. Observe that all the sets in $R$ are disjoint and moreover 
\begin{equation*}
\bigcup_{M\in R} M =[0,1]^d. 
\end{equation*} 
Observe that \begin{eqnarray*}
2^{l_i-1}\delta^{1/d} < h_i \leq 1/2  \;\; \Rightarrow  \;\;2^{l_i} \delta^{1/d} <1 &\Rightarrow & cu^{-1}\delta^{-1/d}2^{-l_i} > 1 \\ & \Rightarrow & \lceil cu^{-1}\delta^{-1/d}2^{-l_i} \rceil \leq 2cu^{-1}\delta^{-1/d}2^{-l_i}.
\end{eqnarray*}
Hence we can easily see that $$|R|=\Pi_{i=1}^d \lceil cu^{-1}\delta^{-1/d}2^{-l_i} \rceil  \leq 2^d c^d u^{-d}\delta^{-1} 2^{-\sum_{i=1}^d l_i}\leq  2^{2d+1} c^d u^{-d}\delta^{-1}. $$ Here the last inequality follows from the fact that $\sum_{i=1}^d l_i \geq -(d+1)$.

Let us define the following set: 
\begin{eqnarray*}
R_2:=\Big\{(M_{\underset{\sim}{i}},M_{\underset{\sim}{i}^\prime})\in R \times R: \exists \; (t,h)\in \mathscr{F}^{\prime} \mbox{ such that } t-h \in M_{\underset{\sim}{i}} \mbox{ and } t+h \in M_{\underset{\sim}{i}^\prime}  \Big\}.
\end{eqnarray*}
%By the very definition of $R_2$ and~\eqref{full}, we have $|R_2| \leq |\mathscr{F}^{\prime}|$. 
Note that if $(t,h)\in \mathscr{F}^{\prime}$ then $h_k \leq 2^{l_k}\delta^{1/d}$ for all $k=1,\ldots,d$. This implies that if $(M_{\underset{\sim}{i}},M_{\underset{\sim}{i}^\prime}) \in R_2$, where $\underset{\sim}{i} = (i_1,\ldots, i_d)$ and $\underset{\sim}{i}^\prime = (i^\prime_1, \ldots, i^\prime_d)$, then 
\begin{equation}\label{eq:i'_k-i_k}
(i^\prime_k- i_k) \leq (1+  2c u^{-1} ), \qquad \mbox{ for all } \;k=1,\ldots,d,
\end{equation} 
as (i) $(i^\prime_k - 1)  u\delta^{{1}/{d}}2^{l_k} c^{-1} \le t_k + h_k$, and (ii) $i_k  u\delta^{\frac{1}{d}}2^{l_k} c^{-1} \ge t_k - h_k$. Thus for each hyperrectangle $M_{\underset{\sim}{i}} \in R$ the number of hyperrectangles $M_{\underset{\sim}{i}^\prime} \in R$ such that $(M_{\underset{\sim}{i}},M_{\underset{\sim}{i}^\prime}) \in R_2$  is less than or equal to $(1+2c u^{-1})^d \leq 4^d c^d u^{-d}$. Hence we have $$|R_2|\leq |R| \times 4^d c^d u^{-d} \leq  2^{4d+1} c^{2d} u^{-2d}\delta^{-1}  \leq d^{2d} 2^{4d^2+4d+1} u^{-2d}\delta^{-1}.$$

Thus, our proof will be complete if we can show that $|R_2|=|\mathscr{F}^\prime|$. From the definition of $R_2$ and the fact that elements in $R$ are disjoint it is easy to observe that $|R_2|\leq |\mathscr{F}^\prime|$. \\

Therefore, the only thing left to show is that $|\mathscr{F}^\prime|\leq |R_2|$. Let us assume the contrary, i.e., $|R_2| < |\mathscr{F}^\prime|$. This implies that there exist two elements $(t,h) $ and $(t^\prime,h^\prime) \in \mathscr{F}^{\prime}$ and $(M_{\underset{\sim}{i}},M_{\underset{\sim}{i}^\prime}) \in R_2 $ such that both $t-h$ and $t^\prime-h^\prime$ belong to $M_{\underset{\sim}{i}}$ and, also, $t+h$ and $t^\prime+h^\prime$ belong to $M_{\underset{\sim}{i}^\prime}$. Let us first define the following two hyperrectangles: 
$$B_1:=\Pi_{k=1}^d (i_k-1,i^\prime_k] \times c^{-1}u\delta^{1/d}2^{l_k} \qquad \mbox{ and } \qquad B_2:=\Pi_{k=1}^d (i_k,i^\prime_k-1] \times c^{-1}u\delta^{1/d}2^{l_k}.$$ Our goal is to show that 
\begin{equation}\label{eq:B_infty} 
B_\infty(t,h) \bigtriangleup B_\infty(t^\prime,h^\prime)\subseteq B_1 \setminus B_2
\end{equation} 
which is implied by the following two assertions: 
\begin{enumerate}
\item[(1)] $B_\infty(t,h) \cup  B_\infty(t^\prime,h^\prime) \subseteq  B_1$ and
\item [(2)] $B_2 \subseteq B_\infty(t,h)\cap B_\infty(t^\prime,h^\prime)$.

\end{enumerate}
See the figure below for a visual illustration of~\eqref{eq:B_infty} when $d=2$. 
\begin{figure}
\includegraphics[scale=.48]{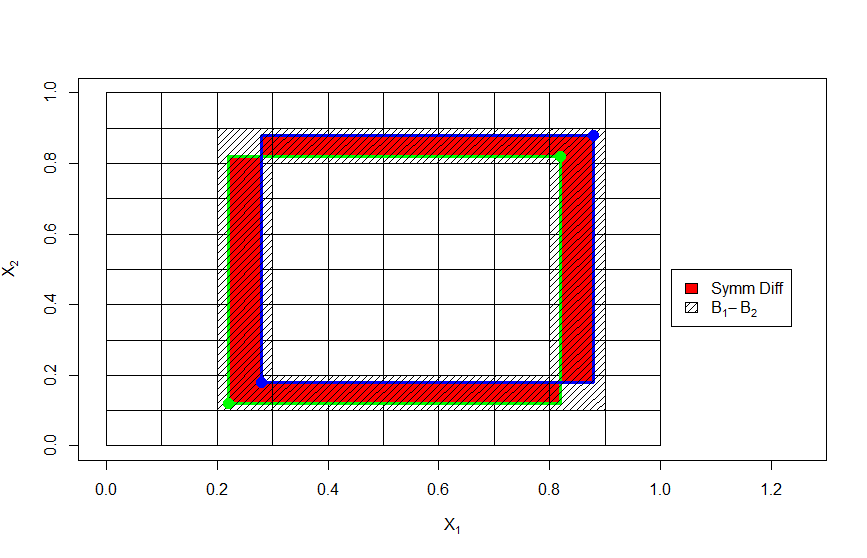}\caption{The figure shows how the symmetric difference of the hyperrectangles $B_\infty(t,h)$ (denoted by the green border) and $B_\infty(t^\prime,h^\prime)$ (denoted by the blue border) is contained in the set $B_1 \setminus B_2$ (denoted by the shaded region). }
\end{figure}
Now, as $t-h\in M_{\underset{\sim}{i}}$, this implies $t_k-h_k \ge (i_k-1)c^{-1}u\delta^{1/d}2^{l_k}$, for all $k=1,\ldots,d$. Also $t+h \in M_{\underset{\sim}{i}^\prime}$ implies that $t_k+h_k \leq i^\prime_{k}c^{-1}u\delta^{1/d}2^{l_k} $, for all $k=1,\ldots,d.$ Therefore, $B_\infty(t,h) = \Pi_{i=1}^d(t_i-h_i,t_i+h_i) \subseteq  B_1 $. A similar argument shows that $B_\infty(t^\prime,h^\prime) \subseteq  B_1$. Hence assertion $(1)$ above holds. \newline

Now as $t-h\in M_{\underset{\sim}{i}}$, we have $t_k-h_k \leq i_kc^{-1}u\delta^{1/d}2^{l_k} $, for all $k=1,\ldots,d.$ Also $t+h \in M_{\underset{\sim}{i}^\prime}$ implies that $t_k+h_k \ge (i^\prime_{k}-1)c^{-1}u\delta^{1/d}2^{l_k} $, for all $k=1,\ldots,d$. Hence we have $B_2 \subseteq B_{\infty}(t,h)$. A similar argument shows that $B_2 \subseteq  B_\infty(t^\prime,h^\prime)$. Therefore, assertion $(2)$ is also satisfied. Now let us define the following set \begin{eqnarray*}I:=\big\{{\underset{\sim}{j}}=(j_1,\ldots,j_d) \in \mathbb{N}^d & : & j_k \in (i_k-1,i^\prime_k], \mbox{ for all } k=1,\ldots,d, \\
& & \exists \: l\in\{1,\ldots,d\} \mbox{ such that } j_l= i_l \mbox{ or } i_l^\prime  \big\}.\end{eqnarray*}
Clearly, using~\eqref{eq:i'_k-i_k}, $$ |I| \leq 2d(2+2cu^{-1})^{d-1}.$$
%(here we have used the fact that $(i_k^\prime - i_k)\leq (1+2cu^{-1})$ for all $k=1,\ldots,d$). 
Also see that  $w = (w_1,\ldots, w_d) \in B_1\setminus B_2$ if and only if
\begin{enumerate}
\item[(1)] for every $k=1,\ldots,d$, we have $w_k \in \big(i_k-1, i_k^{\prime}\big] \times c^{-1}u\delta^{1/d}2^{l_k}$ (this is true as $w\in B_1$), 
\item[(2)]  there exists $l \in \{1,2,\ldots,d\}$ such that either $w_l \in \big(i_l-1, i_l\big] \times c^{-1}u\delta^{1/d}2^{l_l}$  or $w_l \in \big(i_l^\prime -1, i_l^\prime \big] \times c^{-1}u\delta^{1/d}2^{l_l}$ (this is true as $w\not \in B_2$ implies that there exist $l$ such that $w_l \not\in (i_l,i_l^\prime-1]\times c^{-1}u\delta^{1/d}2^{l_l}$ and $w \in B_1$ implies that $w_l \in (i_l-1,i_l^\prime] \times c^{-1}u\delta^{1/d}2^{l_l} $).
\end{enumerate} 
Therefore, we see that $$B_1\setminus B_2= \bigcup_{{\underset{\sim}{j}} \in I} M_{{\underset{\sim}{j}}}.$$ Also, note that, $|M_{\underset{\sim}{j}}| \le u^d\delta c^{-d}2^{\sum_{i=1}^d l_i} \leq u^d\delta c^{-d}$  for all ${\underset{\sim}{j}}$. 
Therefore, using~\eqref{eq:B_infty} and the fact that $c=d4^d$, we easily see that $$\rho^2((t,h),(t^\prime,h^\prime)) \leq |B_1\setminus B_2| \leq  2d (2+2cu^{-1})^{d-1} \times \frac{u^d\delta}{c^d} \le 2^d d (1 + c^{-1})^{d-1} u \delta c^{-1} < u\delta $$
which contradicts \eqref{packingcond}. This proves that two elements of $\mathscr{F}^{\prime}$ cannot correspond to the same pair of hyperrectangles  $(M_{\underset{\sim}{i}},M_{\underset{\sim}{i}^\prime}) \in R_2$. Hence we have proved \eqref{eq3}.\newline

\noindent{\bf Step 2: }
In this part of the proof we show that
\begin{equation}\label{eq4}
N\left((u\delta)^{1/2},\mathscr{F}(\delta)\right)   \leq K u^{-2d}\delta^{-1} (\log(e/\delta))^{d-1}.
\end{equation}
Let us define the set {\small $$S:=\Big\{(l_1,\ldots,\l_d) \in\mathbb{Z}^d:-(d+1) < \sum_{k=1}^d l_k \leq 0  \mbox{ and } l_k \leq \frac{1}{d} \log_2(1/\delta) \mbox{ for all } k=1,\ldots ,d  \Big\}.$$}
Now it can be easily seen that $l:=(l_1,\cdots,\l_d) \in S$ implies $l_k \geq -(d+1) - (d-1)(1/d)\log_2(1/\delta)$, for all $k=1,\ldots,d$. This shows that each $l_k$ can only take at most $(d+2)+ \log_2(1/\delta) \leq (d+2)+\log(1/\delta) \log_2(e) \leq d+2(\log(e/\delta))$ many values. This shows that $$|S| \leq (d+1) (d+2\log(e/\delta))^{d-1} \leq (d+2)^d(\log(e/\delta))^{d-1}.$$
Note that the power of $(d+2\log(e/\delta))$ in the above display is $d-1$ because if we fix the values of $l_1,l_2,\ldots, l_{d-1}$ then $l_d$ can only take at most  $(d+1)$ values such that $(l_1,l_2,\ldots l_{d}) \in S$ (as $\sum_{k=1}^d l_k$ can take at most $d+1$ distinct values). Also note that 
\begin{equation*} 
\mathscr{F}(\delta)\subseteq \bigcup_{l \in S} \mathscr{F}_{\delta,l}.
\end{equation*}
The above representation of $\mathscr{F}(\delta)$ along  with  the trivial fact that $N(\epsilon,\bigcup_{i=1}^n A_i)\leq \sum_{i=1}^n N(\epsilon,A_i)$ gives us \eqref{eq4}.\\[.4cm]
{\bf Step 3: }
In this step we will complete the proof of Lemma~\ref{lem1}. We want control the $\sqrt{u\delta}$-packing number of the set $\{(t,h)\in \mathscr{F}: \sigma^2(t,h) \leq \delta\} $ which can be decomposed in the following way: for $u \in (0,1]$, 
$$\{(t,h)\in \mathscr{F}: \sigma^2(t,h) \leq \delta\} = \left(\bigcup_{l=0}^{\lfloor 1+\log_2(1/u)\rfloor} \mathscr{F}(\delta2^{-l})\right) \cup \{a \in \mathscr{F}: \sigma^2(a) \leq u\delta/2\}.$$
Now we can control the $\sqrt{u\delta}$-packing number of each of the above sets. First observe that
$$N((u\delta)^{1/2},\{(t,h) \in \mathscr{F}: \sigma^2(t,h) \leq u\delta/2\})=1. $$
Also, for any $u\in(0,2)$ and $\delta\in(0,1]$ we have 
\begin{equation}\label{eq4.2}
N((u\delta)^{1/2},\mathscr{F}(\delta)) \leq N((u\delta/2)^{1/2},\mathscr{F}(\delta))\leq K u^{-2d}\delta^{-1} (\log(e/\delta))^{d-1}
\end{equation}
for some constant $K$. Putting $\delta\leftarrow\delta/2^l$ and $u\leftarrow2^lu$ for $0\le l\leq \lfloor 1+\log_2(1/u)\rfloor$ in \eqref{eq4.2} we get $$N((u\delta)^{1/2},\mathscr{F}(\delta2^{-l}) ) \leq K 2^{-(2d-1)l} u^{-2d} \delta^{-1}(\log(e/\delta))^{d-1}.$$ Now from the trivial fact that $N(\epsilon,\bigcup_{i=1}^m A_i)\leq \sum_{i=1}^m N(\epsilon,A_i)$  we get 
\begin{eqnarray*}
&& N\left(\sqrt{u\delta},\{(t,h)\in \mathscr{F}: \sigma^2(t,h) \leq \delta\}\right) \\
& \leq & \sum_{l=0}^{\lfloor 1+\log_2(1/u)\rfloor}N\left(\sqrt{u\delta},\mathscr{F}(\delta2^{-l})\right) + N\left(\sqrt{u\delta},\{(t,h) \in \mathscr{F}: \sigma^2(t,h) \leq u\delta/2\}\right)\\
%& \leq & 1+\sum_{l=0}^{\lfloor 1+\log_2(1/u)\rfloor} K 2^{-(2d-1)l}u^{-2d}\delta^{-1} (\log(e/\delta))^{d-1}\\
& \leq & 1+ Ku^{-2d}\delta^{-1} (\log(e/\delta))^{d-1} \sum_{l=0}^\infty 2^{-(2d-1)l}\\
& \leq & 1+ 2Ku^{-2d}\delta^{-1} (\log(e/\delta))^{d-1}\\
& \leq &(2K+1)u^{-2d}\delta^{-1} (\log(e/\delta))^{d-1},
\end{eqnarray*}
which proves Lemma~\ref{lem1}.  \qed \newline
\subsection{Proof of Theorem~\ref{Thm 1}}\label{sec:Thm 1}
We use Theorem~\ref{Thm0} to prove  Theorem~\ref{Thm 1}. 
%\noindent \textbf{Proof of Theorem 2.1.}\\ 
Let us recall the definitions of $\mathscr{F}, \sigma$ and $\rho$ as introduced just before Lemma~\ref{lem1}. Without loss of generality we assume that $\norm{\psi}=1$. For $h \in (0,1/2]^d$, let us define the stochastic process
$$X(t,h) :=2^{d/2}(h_1h_2 \ldots h_d)^{1/2} \hat{\Psi}(t,h)=2^{d/2} \int \psi_{t,h}(x) dW(x), \qquad t \in A_h,$$ where $W(\cdot)$ is the standard Brownian sheet on $[0,1]^d$. This defines a centered Gaussian process with $\Var\big(X(t,h)\big)=\sigma^2(t,h)$. Also by a standard calculation on the variance we have $\Var\big(X(t,h)-X(t^\prime,h^\prime)\big) \leq 2^d TV^2(\psi)\rho^2((t,h),(t^\prime,h^\prime))$. 
  As $X(t,h)$ and $X(t,h)-X(t^\prime,h^\prime)$  have normal distributions this shows that conditions (a) and (b) of Theorem \ref{Thm0} are satisfied. Condition (c) is also satisfied because of Lemma \ref{lem1}.  Thus, by an application of Theorem \ref{Thm0} we have 
$$\sup_{0<h\leq1/2}\sup_{t\in A_{h}} \frac{\hat{\Psi}^2(t,h)- 2\log(1/2^dh_1h_2...h_d)}{\log\log(e^e/2^dh_1h_2...h_d)} < \infty.$$

For notational simplicity, let us define $\kappa_1 :=2\log(1/\sigma^2(t,h))$~and~$\kappa_2 :=2\sqrt{2} S\log\log({e^e/\sigma^2(t,h)})$. Therefore,
\begin{eqnarray*}
& & {\small \p \left(|\hat{\Psi}(t,h)| \leq \sqrt{2\log\left(\frac{1}{\sigma^2(t,h)}\right)}+  S \Bigg(\frac{ \log\log(e^e/\sigma^2(t,h))}{\log^{\frac{1}{2}}(1/\sigma^2(t,h)) } \Bigg) \;\; \forall \; (t,h) \in \mathscr{F} \right)}\\
&=&\p\left(|\hat{\Psi}(t,h)| \leq \kappa_1^{1/2}+ \kappa_1^{-1/2} \kappa_2/2 \quad \forall \; (t,h) \in \mathscr{F} \right)\\
&=& \p \left(\hat{\Psi}(t,h)^2 \leq \Big(\kappa_1^{1/2}+ \kappa_1^{-1/2} \kappa_2/2\Big)^2 \: \forall (t,h) \in \mathscr{F} \right)\\
&\geq & \p \left(\hat{\Psi}(t,h)^2 \leq  \kappa_1+\kappa_2 \quad \forall (t,h) \in \mathscr{F} \right)\\
&=& \p \left(\sup_{t,h\in \mathscr{F}} \frac{\hat{\Psi}^2(t,h)- 2\log(1/2^dh_1h_2...h_d)}{\log\log(e^e/2^dh_1h_2...h_d)} <  2\sqrt{2}S \right) \; \to \;  1 \quad \mbox{  as    }S\to \infty. \qed
\end{eqnarray*}

\appendix
\section{Proofs of other results}\label{proofs of result2}
\subsection{Proof of Proposition \ref{v<1}}\label{sec:Prop_V}
The proof of this  result follows from the following result. Suppose that $Z_1,Z_2,\ldots ,Z_n$ are i.i.d.~standard normal random variables. Then, we know that $$\frac{\max_{1\leq i\leq n} Z_i}{\sqrt{2\log n}} \rightarrow 1 \quad \mbox{ a.s}. $$
Let $F_n$ be the distribution function of $\max_{1\leq i\leq n} Z_i/\sqrt{2\log n}$, i.e.,  $F_n(x):= \p({\max_{1\leq i\leq n} Z_i} \leq x \sqrt{2\log n})$, for $x \in \R$. Therefore, for every $x<1$, we have $F_n(x)\to 0.$ We want to show that $$\sup_{(t,h)\in \mathscr{F}} |\hat{\Psi}(t,h)|-\Gamma_V(2^dh_1\ldots h_d)=\infty \quad \mbox{a.s.}$$ Hence it is enough to show that for every $s\in \R$ we have $\p(\sup_{(t,h)\in \mathscr{F}} |\hat{\Psi}(t,h)|-\Gamma_V(2^dh_1\ldots h_d)<s)=0$. Fix $m \in \N$. Now, 
\begin{eqnarray*}
& &\p\left(\sup_{(t,h)\in \mathscr{F}} |\hat{\Psi}(t,h)|-\Gamma_V(2^dh_1\ldots h_d)<s\right)\\
& \leq & \p\left(\sup_{t\in A_{\left(\frac{1}{2m},\ldots,\frac{1}{2m}\right)}} \left|\hat{\Psi}\left(t,\left(\frac{1}{2m},\ldots,\frac{1}{2m}\right)\right)\right|-\Gamma_V(m^{-d})<s\right)\\
&\leq & \p\left(\sup_{t\in A_m^\star} |\hat{\Psi}(t,({2m})^{-1})|-\Gamma_V(m^{-d})<s\right) 
\end{eqnarray*}
where $A_m^\star :=\{(t_1,\ldots,t_d): t_i=k_i/2m \mbox{ for some odd integer  } k_i <2m, \; \mbox{for all }\; i = 1,\ldots, d\}$. Thus, the last term in the above display can be further upper bounded by 
\begin{eqnarray*}
\p\left(\sup_{t\in A_m^\star} \frac{\hat{\Psi}(t,(2m)^{-1})}{\sqrt{2\log (m^d)}}- \sqrt{V} < \frac{s}{\sqrt{2\log (m^d)}}   \right) = F_{m^d}(\sqrt{V} + s/\sqrt{2\log (m^d)} ),
\end{eqnarray*} 
where we have used the fact that now we are dealing with $m^d$ i.i.d.~standard normal random variables. Now, for every $s>0$, choose $m$ such that $\sqrt{V} + s/\sqrt{2\log (m^d)}< 1-\epsilon$, for some fixed $\epsilon>0$. Hence, $F_{m^d}(\sqrt{V} + s/\sqrt{2\log (m^d)} )\leq F_{m^d}(1-\epsilon)$, if $m$ is large enough. As this is true for all large $m$, taking $m \to \infty$ gives us the desired result. \qed
\subsection{Solution to \eqref{minimize}}\label{lem:Ele}
Let $\psi \in \mathbb{H}_{\beta,1}$ such that $\psi(0) \geq 1.$ Hence by the property of $\mathbb{H}_{\beta,1}$ we have $$|\psi(x)-\psi(0)| \leq \norm{x}^\beta, \qquad \mbox{for all } x \in \R^d,$$ which implies $\psi(x) \geq 1-\norm{x}^\beta$. Hence, on the set $\norm{x} \leq 1$, we have $\psi(x) \geq 1-\norm{x}^\beta \geq 0.$ Therefore, we have $$\int_{\norm{x} \leq 1} \psi^2(x) dx \geq \int_{\norm{x}\leq 1} (1-\norm{x}^\beta)^2 dx  \quad \Rightarrow \quad \norm{\psi} \geq  \norm{\psi_\beta},$$
%And also obviously we have $$\int_{\norm{x} \geq 1} \psi^2(x) dx \geq 0.$$ The above two displays gives us that $\norm{\psi} \geq  \norm{\psi_\beta}$ 
where $\psi_\beta(x)=(1-\norm{x}^\beta) \mathbb{I}({\norm{x} \leq 1}).$
Hence the only thing left to prove is that $\psi_\beta \in \mathbb{H}_{\beta,1}$.
Suppose that $x,y \in \mathbb{R}^d$ such that $1\geq \norm{x} \geq \norm{y}$. Then \begin{equation*}
\begin{split}
0\leq \psi_\beta(y)-\psi_\beta(x) =\norm{x}^\beta -\norm{y}^\beta  \leq (\norm{x} -\norm{y})^\beta \leq \norm{x-y}^\beta.
\end{split}
\end{equation*}
Here the the third inequality follows from the fact that when $\beta\leq 1$ the function $u \mapsto u^\beta$ is a $\beta$-H\"older continuous function; the last inequality follows from the triangle inequality.
If $x,y \in \mathbb{R}^d$ such that $ \norm{x}\geq 1 \geq \norm{y}$ then we have 
\begin{equation*}
\begin{split}
0\leq \psi_\beta(y)-\psi_\beta(x) =1 -\norm{y}^\beta  \leq (1 -\norm{y})^\beta \leq (\norm{x} -\norm{y})^\beta \leq \norm{x-y}^\beta.
\end{split}
\end{equation*}
If $x,y \in \R^d$ is such that $\norm{x}\geq \norm{y}\geq 1$ then the assertion is trivial.
Hence we have proved that $\psi_\beta$ minimizes~\eqref{minimize}. \qed

\subsection{Proofs of Theorems~\ref{multopt} and~\ref{boxthm}}\label{multiopt}
The proofs of~Theorems~\ref{multopt} and \ref{boxthm} depend on the following lemma (stated and proved in~\citet[Lemma 6.2]{DC01}).
\begin{lem}\label{normal lemma}
Let $Z_1,Z_2,\ldots$ be a sequence of independent standard normal variables. If $w_m := (1-\epsilon_m)\sqrt{2\log m}$ with $\lim_{m \to \infty} \epsilon_m =0$ and $\lim_{m \to \infty} \epsilon_m \sqrt{\log m}=\infty$, then we have $$\lim_{m \to \infty} \E\left|\frac{1}{m}\sum_{i=1}^m \exp\left(w_mZ_i-\frac{w_m^2}{2}\right) - 1\right|=0.
$$
\end{lem}
%The proof of the above lemma can be found in \citet{DC01}. \newline

\subsubsection{Proof of Theorem \ref{multopt}}
\noindent \textit{Proof of part} (a).  For any bandwidth $h = (h_1,\ldots, h_d)  \in (0,1/2]^d$ and $t=(t_1,\ldots,t_d)\in A_h$, let us define the function $g_t:[0,1]^d \to \R$ as $$g_t(x) := L \min(h)^\beta \psi^{(\beta)}_{t,h}(x), \quad \mbox{for } x \in [0,1]^d,$$ where $\min(h):=\min\{h_1,h_2, \ldots,h_d\}$ and $\psi_{t,h}^{(\beta)}(x_1,\ldots,x_d)=\psi_\beta((x_1-t_1)/h_1,\ldots,(x_d-t_d)/h_d)$. Elementary calculations show that $g_t \in\mathbb{H}_{\beta,L}$ and $\|g_t\|_{\infty}=L \min(h)^{\beta}$. Now let us define the set $$S :=\big\{t \in A_h: t_i= k_i h_i \mbox{ for some odd integer } k_i, i = 1,\ldots, d\big\}.$$ 

Let $\phi_n$ be an arbitrary test for~\eqref{test} with level $\alpha$. Then,
\begin{eqnarray}
&\inf_{g \in \mathbb{H}_{\beta,L}: \norm{g}_\infty= L\min(h)^\beta} \E_g [\phi_n{(Y)}] - \alpha & \; \leq  \min_{g_t: t\in S} \E_{g_t} [\phi_n(Y)]- \E_0 [\phi_n(Y)] \nonumber\\
&&\leq \;  |S|^{-1} \sum_{t \in S} \E_{g_t} [\phi_n(Y)]-\E_0[\phi_n(Y)]\nonumber\\
&&\leq \; \E_0\Bigg[\Big(|S|^{-1}\sum_{t \in S } \frac{dP_{g_t}}{dP_{0}}(Y) -1\Big)\phi_n(Y)\Bigg]\nonumber\\
&&\leq \; \E_0 \Big| |S|^{-1}\sum_{t \in S } \frac{dP_{g_t}}{dP_{0}}(Y) -1\Big| \label{eq:imp}.
\end{eqnarray}
Here $P_0$ denotes the measure of the process $Y$ under the null hypothesis $f=0$ and $P_{g_t}$ denotes the measure of $Y$ under the alternative $f=g_t$. Also  for $g \in \mathbb{H}_{\beta,L}$, $\frac{dP_{g}}{dP_{0}}$ denotes the Radon-Nikodym derivative of the measure $P_{g}$ with respect to the measure $P_0$. By Cameron-Martin-Girsanov's Theorem (see \citet[Chapter 3]{Protterbook} for more details about absolute continuous measures and Radon-Nikodym derivatives) we get that  $$\log\left(\frac{dP_{g}}{dP_{0}}(Y)\right)= \sqrt{n} \int g dW - \frac{n}{2} \norm{g}^2.$$ For $g_t(\cdot) =L\min(h)^\beta \psi^{(\beta)}_{t,h}(\cdot)$, $\sqrt{n} \int g_t dW=\sqrt{n}L \norm{\psi_\beta} \min(h)^\beta \sqrt{\Pi_{i=1}^d h_i } \hat{\Psi}(t,h).$ Observe that $\{Z_t \equiv \hat{\Psi}(t,h)\}_{t \in S}$ are i.i.d.~standard normals. Let $$w_n:=\sqrt{n}L \norm{\psi_\beta} \min(h)^\beta \sqrt{\Pi_{i=1}^d h_i}.$$ Then $\Gamma_t=\exp(w_nZ_t-\frac{w_n^2}{2})$ and we can write $\frac{dP_{g_t}}{dP_0}(Y)-1=\Gamma_t-1.$

Hence we have $\E_0 \left| |S|^{-1} \sum_{t \in S } \frac{dP_{g_t}}{dP_{0}}(Y) -1\right| = \E_0  \left| |S|^{-1} \sum_{t \in S }  \Gamma_t -1  \right| $. According to Lemma~\ref{normal lemma} the above term will go to zero if $|S|\to \infty$ and the corresponding $w_n$'s satisfy: 
$$\left(1-\frac{w_n}{\sqrt{2\log |S| }}\right)  \to 0  \qquad \mbox{and} \qquad \sqrt{\log |S| } \left(1-\frac{w_n}{\sqrt{2\log |S| }}\right) \to \infty.$$

Now let us pick $$h_1= \ldots =h_d= L^{-\frac{2}{2\beta+d}}((1-\epsilon_n)\rho_n)^{1/\beta} \left( \norm{\psi_\beta}^2(2\beta+d)/2d\right)^{-1/(2\beta+d)} =: \tilde h.$$
Then, 
\begin{eqnarray}\label{6.18} w_n &=&\sqrt{n}L\norm{\psi_\beta} L^{-1}\left((1-\epsilon_n)\rho_n\right)^\frac{2\beta+d}{2\beta} \left( \norm{\psi_\beta}^2(2\beta+d)/2d\right)^{-1/2} \nonumber\\
&=& \sqrt{n} (1-\epsilon_n)^{1+d/2\beta}\sqrt{\frac{\log n}{n}}\sqrt{(2d/(2\beta+d))}\nonumber \\
& = & \sqrt{(2d/(2\beta+d))} (1-\epsilon_n)^{1+d/2\beta} \sqrt{\log n}.
\end{eqnarray}
Also, as $n\to \infty$, $|S|/(\Pi_{i=1}^d(1/ h_i)) \to 2^{-d}$. Therefore, for a suitable constant $K$, 
\begin{eqnarray}\label{6.19}
\log |S|/\log n&=&(- d \log \tilde h -d\log 2 +o(1))/\log n \nonumber\\
&= &[K+o(1) - ({d}/{\beta})\log\left((1-\epsilon_n)\rho_n\right)]/ \log n  \nonumber\\
& = & \left(K+o(1)-\frac{d}{\beta}\log(1-\epsilon_n)+ \frac{d}{2\beta+d}\log \left(\frac{n}{\log n}\right)\right)/\log n \nonumber \\
& \to & \frac{d}{2\beta+d} \mbox{    as      } n\to \infty. \end{eqnarray} 
Also notice that for all large $n$,
$\log |S|/\left(\frac{d}{2\beta+d}\log n\right)<1$. Combining~\eqref{6.18} and~\eqref{6.19}, we get  $$\frac{w_n}{\sqrt{2\log |S|}} = \frac{w_n}{\sqrt{\log n}} \frac{\sqrt{\log n}}{\sqrt{2\log |S|}} \to 1\quad \mbox{as } \; n \to \infty.$$ Similarly, for suitable constants $K, K'>0$, 
\begin{eqnarray*} \sqrt{\log |S| } \left(1-\frac{w_n}{\sqrt{2\log |S| }}\right)&\geq &\sqrt{K}\sqrt{\log n}\left(1-(1-\epsilon_n)^{1+d/2\beta} + o(1)\right)\\
&\geq &\sqrt{K'}\sqrt{\log n}\left(\epsilon_n+o(1)\right) \to \infty \quad \mbox{as } \; n \to \infty,
\end{eqnarray*}
as the o(1) term above is positive when $n$ is large. This proves part (a) of Theorem~\ref{multopt} by noting that  $L\min(h)^\beta=(1-\epsilon_n)c_*\rho_n$. \newline

\noindent \textit{Proof of part} (b). 
Let $\delta \equiv \delta_n:=c_*\rho_n$ and $h_{i,n}=(\delta/L)^{1/\beta} =: \tilde h_n$ for all $i=1,2,\ldots,d$. For notational simplicity, in the following we drop the subscript $n$. As the term $D(2^dh_1\ldots h_d)$ is bounded from above, for any $t \in J$, the probability of rejecting the null hypothesis, $\p_g(T_\beta(Y) > \kappa_\alpha)$, is bounded from below by, for some constant $K >0$, 
\begin{eqnarray}\label{eq:powerTb}
& \quad   & \quad \: \p_g\left(|\hat{\Psi}(t,h)| > \Gamma(2^d \tilde h^d) + K \right) \nonumber \\
& \quad & = \, \p_0 \left(\left| \hat{\Psi}(t,h) + \sqrt{\frac{n}{\tilde h^d}} \norm{\psi_\beta}^{-1} \langle g, \psi^{(\beta)}_{t,h}\rangle  \right|> \Gamma(2^d \tilde h^d) + K\right) \nonumber \\
& \quad & \geq   \p_0\left(-\mbox{sign}(\langle g, \psi^{(\beta)}_{t,h}\rangle)\hat{\Psi}(t,h) < \sqrt{\frac{n}{\tilde h^d}} \frac{ |\langle g, \psi^{(\beta)}_{t,h}\rangle |}{\norm{\psi_\beta}} - K - \Gamma(2^d \tilde h^d)\right) \nonumber\\
 & \quad &= \,\Phi\left(\sqrt{\frac{n}{\tilde h^d}} \norm{\psi_\beta}^{-1} |\langle g, \psi^{(\beta)}_{t,h}\rangle | - K - \Gamma(2^d \tilde h^d)\right)
\end{eqnarray}
where $\Phi$ is the standard normal distribution function. Hence, to prove our claim it suffices to show that $$(1+\epsilon_n) \max_{t \in J} \sqrt{\frac{n}{\tilde h^d}} \norm{\psi_\beta}^{-1} |\langle g, \psi^{(\beta)}_{t,h}\rangle | - \Gamma(2^d \tilde h^d) \to \infty$$ uniformly for all $g \in \mathbb{H}_{\beta,L}$ such that $\norm{g}_{J,\infty} \ge \delta$.
Note that $A_{h}=J.$

Let $g$ be any such function, and let $t \in J$ be such that $|g(t)| \ge \delta$. Let us assume that  $g(t)\ge \delta$; the other case where $g(t)\le -\delta$ can be handled similarly by looking at $-g$. By construction of $\psi_\beta$ we have $\delta \psi_{t,h}^{(\beta)} \in \mathbb{H}_{\beta,L}$. Also note that as $\psi_\beta$ minimizes $\norm{\psi}$ in the set $\{\psi \in \mathbb{H}_{\beta,1}: \psi(0)\geq 1\}$, $\delta\psi_{t,h}^{(\beta)}$ minimizes $\norm{\psi}$ in the set $\{\psi \in \mathbb{H}_{\beta,L}: \psi(t)\geq \delta\}$. Note that both $g$ and $\delta\psi_{t,h}^{(\beta)}$ belong to the closed convex set $\{\psi \in \mathbb{H}_{\beta,L}: \psi(t)\geq \delta\}$. As $\delta\psi_{t,h}^{(\beta)}$ is the projection of the zero function onto the above closed convex set, we have  $$|\langle \psi_{t,h}^{(\beta)},g\rangle| = \delta^{-1} |\langle \delta\psi_{t,h}^{(\beta)},g\rangle| \geq  \delta^{-1} \|\delta \psi_{t,h}^{(\beta)}\|^2= \delta \norm{\psi_\beta}^2 \tilde h^d.$$ Thus, 
\begin{equation*}
\begin{split}
& (1+\epsilon_n) \max_{t \in J} \sqrt{\frac{n}{\tilde h^d}} \norm{\psi_\beta}^{-1} |\langle g, \psi^{(\beta)}_{t,h}\rangle | - \Gamma(2^d \tilde h^d)\\
& \geq \, (1+\epsilon_n)\norm{\psi_\beta}\delta \sqrt{n \tilde h^d} -\Gamma(2^d \tilde h^d)\\
& =\, (1+\epsilon_n)\norm{\psi_\beta} c_*\rho_n \sqrt{n} (c_*\rho_n)^{d/2\beta}L^{-d/2\beta} - \Gamma(2^d \tilde h^d)\\
& = \, (1+\epsilon_n)\sqrt{\left(\frac{2d}{2\beta+d}\right)\log n}- \sqrt{K+\left(\frac{2d}{2\beta+d}\right)\log\left(\frac{n}{\log n}\right)}\\
& \geq\, \epsilon_n (2d/(2\beta+d))^{1/2}(\log n)^{1/2} + o(1) \to \infty.
\end{split} 
\end{equation*}
This proves part (b) of  Theorem \ref{multopt}.\qed%\vspace{.3cm}

\subsubsection{Proof of Proposition~\ref{prop:triangle}}\label{proof:proptriangle} 
Let $h :=(\tilde{h},\ldots,\tilde{h}) \in \R^d$, where $\tilde{h}= (M\rho_n/L)^{1/\beta}$, for $M$ as defined in the statement of the proposition. 
%\begin{equation}\label{inequ:M} M >\left(\frac{2dL^{d/\beta}\norm{\psi_1}^2}{(2\beta+d)\langle \psi_1,\psi_\beta \rangle^2}\right)^{\frac{\beta}{2\beta+d}}.
%\end{equation}
By the same argument as in~\eqref{eq:powerTb} we have 
$$\p_g(T(Y) > \kappa_\alpha) \geq \,\Phi\left(\sqrt{\frac{n}{\tilde h^d}} \norm{\psi_1}^{-1} |\langle g, \psi^{(1)}_{t,h}\rangle | - K - \Gamma(2^d \tilde h^d)\right).$$
Now we would want to bound $|\langle g, \psi^{(1)}_{t,h}\rangle |$ uniformly for all $g\in\mathbb{H}_{\beta,L}$ such that $\norm{g}_{J_n,\infty} \geq M\rho_n $. Without loss of generality, let us assume that $g(t) \geq M\rho_n $ for some $t \in J_n$ and $g\in \mathbb{H}_{\beta,L}$. Then $$g(x) \geq g(t)-L \norm{x-t}^\beta \geq M\rho_n - L \norm{x-t}^\beta = M\rho_n\left(1 -\norm{\frac{x-t}{\tilde{h}}}^\beta\right).$$ 
This shows that if $\norm{x-t} \leq \tilde{h}$ then $g(x)\geq 0.$
Hence,
\begin{eqnarray*}
\langle g,\psi^{(1)}_{t,h}\rangle &\geq& \int_{\norm{x-t} \leq \tilde{h}} M\rho_n\left(1 -\norm{\frac{x-t}{\tilde{h}}}^\beta\right) \left(1-\norm{\frac{x-t}{\tilde{h}}}\right) dx\\
&=&  M\rho_n\tilde{h}^d\int_{\norm{x}\leq 1} \left(1-\norm{x}\right) \left(1-\norm{x}^\beta \right)dx\\
&=&  M\rho_n\tilde{h}^d\langle\psi_\beta,\psi_1\rangle.
\end{eqnarray*}
Here the last equality follows as $\psi_\beta(x)= (1-\norm{x}^\beta)\mathbb{I}(\norm{x}\leq 1)$. 
Also note that {\small $$\Gamma(2^d\tilde{h}^d)= \sqrt{2d \log\left(\frac{1}{2}\right)+ \frac{2d}{\beta} \log\left(\frac{L}{M}\right) + \frac{2d}{2\beta+d} \log\left( \frac{n}{\log n}\right)} \leq \sqrt{\frac{2d}{2\beta+d} \log n}$$} for large $n$.
Therefore, for large $n$,  \begin{eqnarray*}
& &\sqrt{\frac{n}{\tilde h^d}} \norm{\psi_1}^{-1} \langle g, \psi^{(1)}_{t,h}\rangle  - K - \Gamma(2^d \tilde h^d)\\ & \geq & \sqrt{n\tilde{h}^d} M\rho_n \frac{\langle\psi_\beta,\psi_1\rangle}{\norm{\psi_1}}-K-\sqrt{\frac{2d}{2\beta+d} \log n}\\
& =& -K + \sqrt{\log n} \left(L^{-d/2\beta} M^{\frac{(d+2\beta)}{2\beta}}\frac{\langle\psi_\beta,\psi_1\rangle}{\norm{\psi_1}}- \sqrt{\frac{2d}{2\beta+d}}\right) \to \infty \mbox{   as    } n \to \infty.
\end{eqnarray*}
Here the last equality holds by the choice of $M$, as 
\begin{eqnarray*}
\sqrt{n\tilde{h}^d} M \rho_n \frac{\langle\psi_\beta,\psi_1\rangle}{\norm{\psi_1}} &=& \sqrt{n} M^{\frac{d}{2\beta}} \rho_n^{\frac{d}{2\beta}} L^{-\frac{d}{2\beta}} M \rho_n \frac{\langle\psi_\beta,\psi_1\rangle}{\norm{\psi_1}}\\
%& = & \sqrt{n}  M^{\frac{(d+2\beta)}{2\beta}}L^{-\frac{d}{2\beta}}\rho_n^\frac{d+2\beta}{2\beta}\frac{\langle\psi_\beta,\psi_1\rangle}{\norm{\psi_1}} \\
%& = & \sqrt{n}  M^{\frac{(d+2\beta)}{2\beta}}L^{-\frac{d}{2\beta}}\sqrt{\frac{\log n}{n}}\frac{\langle\psi_\beta,\psi_1\rangle}{\norm{\psi_1}}\\
& = & \sqrt{\log n} \: L^{-d/2\beta} M^{\frac{(d+2\beta)}{2\beta}}\frac{\langle\psi_\beta,\psi_1\rangle}{\norm{\psi_1}}\\
& > & \sqrt{\log n} \: \sqrt{\frac{2d}{2\beta+d}}. \quad\\
\end{eqnarray*}
Hence $\lim_{n \to \infty } \p_g(T(Y) > \kappa_\alpha)=1.$ \qed
\subsubsection{Proof  of Theorem \ref{boxthm}}
\noindent \textit{Proof of part} (a). Let us suppose that $B_n :=B_\infty(t_n,h_n) \subseteq [0,1]^d$ for some $t_n, h_n \in [0,1]^d$. Let us first look at the case when $\liminf_{n \to \infty} |B_n| > 0$. Now assume that the location $B_n$ was known and it was also known that $\mu_n>0$. In such a scenario the best test statistic would be $\hat{\Psi}(t_n,h_n)$ (with kernel $\psi_0$) which follows the normal distribution with mean $0$ and variance $1$, under the null hypothesis. Hence in this case, the UMP test rejects $H_0: \mu_n=0$ if $\hat{\Psi}(t_n,h_n) > z_{1-\alpha}$ where $z_{1-\alpha}$ is the $(1-\alpha)$'th quantile of the standard normal distribution. When $B_n$ is not known then, obviously, the power of any test $\phi_n$ is less than the test described above.
Hence, 
\begin{eqnarray*}
\E_{f_n} [\phi_n(Y)] & \leq & \p_{\mu_n}\left(\hat{\Psi}(t_n,h_n) \geq z_{1-\alpha}\right)=\p_0\left(\hat{\Psi}(t_n,h_n)+ \sqrt{n|B_n|}\mu_n  \geq z_{1-\alpha} \right)\\
& = & 1- \Phi\left(z_{1-\alpha} -\sqrt{n|B_n|} \mu_n   \right)\\
& \not \to &  1 \mbox{  unless  } \mu_n\sqrt{n |B_n|} \to \infty. 
\end{eqnarray*}
A similar argument can be made when $\mu_n < 0$ as well. Hence the power of any test does not go to $1$  unless $|\mu_n|\sqrt{n |B_n|} \to \infty$.

Now suppose that $|\mu_n| \sqrt{n|B_n|}\to \infty $.  Then we will show that $\lim_{n\to \infty} \p_{f_n}(T > \kappa_\alpha)=1$. Without loss of generality assume $\mu_n>0$. Hence,
\begin{eqnarray*}
\p_{f_n}(T> \kappa_\alpha) & \geq & \p_{f_n}\left(\frac{|\hat{\Psi}(t_n,h_n)| - \Gamma(|B_n|)}{D(|B_n|)} > \kappa_\alpha\right)\\
&=& \p_0\left(\left|\hat{\Psi}(t_n,h_n) + \mu_n \sqrt{n|B_n|} \right| - \Gamma(|B_n|)  \geq \kappa_\alpha D(|B_n|)\right)\\
&\geq & \p_0\left(\left|\hat{\Psi}(t_n,h_n) +  \mu_n \sqrt{n|B_n|} \right| \geq K \right) \to 1 \mbox{ as } \mu_n \sqrt{n|B_n|}\to \infty.
\end{eqnarray*}
Here the last inequality follows from the fact that as $ \liminf_n |B_n| > 0$, $\Gamma(|B_n|) + \kappa_\alpha D(|B_n|)$ is bounded from above (say, by $K$) for all large $n$. \newline
%This proves part $(a)$ of our theorem.

\noindent \textit{Proof of part} (b). Now let us look at the case $\lim|B_n| \to 0 $. Let us assume that $|\mu_n|\sqrt{n|B_n|} = (1-\epsilon_n)\sqrt{2\log(1/|B_n|)}$ where $\epsilon_n \to 0$ and $\epsilon_n \sqrt{2\log(1/|B_n|)} \to \infty$. Without loss of generality also assume that $\mu_n>0$. Recall that $B_n  =B_\infty(t_n,h_n)$ for $h_n = (h_{1,n},\ldots, h_{d,n}) \in (0,1/2]^d$. Let us first define the following grid points: {\small $$G_{h_n}:= \left\{t = (t_1,\ldots, t_d) \in [0,1]^d: t_i=(2k_i-1)h_{i,n} \mbox{ for some } k_i \in\N , B_{\infty}(t,h_n) \subseteq [0,1]^d \right\}.$$}
Clearly $|G_{h_n}| \leq 1/|B_n|$. Also, as $n\to \infty$, $|G_{h_n}||B_n| \to 1$. For each $t \in G_{h_n} $  define $f_t := \mu_n \mathbb{I}_{B_\infty(t,h_n)}$. Clearly as $|B_n|=|B_\infty(t,h_n)|$, we have $f_t \in \G_n^{-}$. Let $\phi_n$ be a test of level $\alpha$ for testing~\eqref{boxtest}. Similar arguments as in~\eqref{eq:imp} show that $$\inf_{g \in \G_n^{-}}\E_g \phi_n(Y) - \alpha \leq \E_0\left||G_{h_n}|^{-1} \sum_{t \in G_{h_n}} \frac{dP_{f_t}}{dP_0}(Y) -1
\right|.$$ Now by an argument a similar to that in the proof of Theorem~\ref{multopt}, we have
$$\log\left(\frac{dP_{f_t}}{dP_{0}}(Y)\right)= \sqrt{n} \int f_t dW - n \norm{f_t}^2/2= \mu_n \sqrt{n|B_n|}\hat{\Psi}(t,h_n) - \mu_n^2 n |B_n|/2 .$$ Also note that the collection of random variables in $\{\hat{\Psi}(t,h_n): t \in G_{h_n}\}$ are mutually independent. Now putting $w_n=\mu_n \sqrt{n|B_n|} = (1-\epsilon_n)\sqrt{2 \log(1/|B_n|)} $ and $m=|G_{h_n}|$ we see that $$\E_0\left||G_{h_n}|^{-1} \sum_{t \in G} \frac{dP_{f_t}}{dP_0}(Y) -1
\right| \to 0 $$ if $\epsilon_n \to 0$ and $\epsilon_n \sqrt{\log(1/|B_n|)} \to \infty$, by a direct application of Lemma~\ref{normal lemma}. This proves that
$$ \limsup_{n \to \infty} \inf_{f_n \in \G_n^{-}} \E_{f_n} \phi_n \leq \alpha.$$

Now let us assume that $|\mu_n|\sqrt{n|B_n|} \geq (1+\epsilon_n)\sqrt{2\log(1/|B_n|)}$. Without loss of generality also assume that $\mu_n>0$. A similar argument as in part $(a)$ shows that 
\begin{eqnarray*}
\p_{f_n}(T> \kappa_\alpha) & \geq & \p_{f_n}\left(\frac{|\hat{\Psi}(t_n,h_n)| - \Gamma(|B_n|)}{D(|B_n|)} > \kappa_\alpha\right)\\
&=& \p_0\left(\left|\hat{\Psi}(t_n,h_n) + \mu_n \sqrt{n|B_n|} \right| \geq  \Gamma(|B_n|) + \kappa_\alpha D(|B_n|)\right)\\
& \geq & \p_0\left(\hat{\Psi}(t_n,h_n) \geq  \Gamma(|B_n|) + \kappa_\alpha D(|B_n|) - \mu_n \sqrt{n|B_n|}\right)\\
& \geq & \p_0\left(\hat{\Psi}(t_n,h_n) \geq  -\epsilon_n \sqrt{2\log(1/|B_n|)} + \kappa_\alpha D(|B_n|)\right) \to 1 \mbox{ as } n\to \infty.
\end{eqnarray*}
This completes the proof of Theorem~\ref{boxthm}. \qed

\bibliographystyle{apalike}
\bibliography{Biblio}
\end{document}